\documentclass[a4paper,12pt]{article}%
\usepackage{amsmath}
\usepackage{amsfonts}
\usepackage{amssymb}
\usepackage{graphicx}%
\newtheorem{theorem}{Theorem}

\newtheorem{corollary}{Corollary}
\newtheorem{definition}{Definition}
\newtheorem{example}{Example}
\newtheorem{claim}{Claim}
\newtheorem{lemma}{Lemma}
\newtheorem{proposition}{Proposition}
\newtheorem{remark}{Remark}
\newenvironment{proof}[1][Proof]{\noindent\textbf{#1.} }{\ \rule{0.5em}{0.5em}}
\begin{document}

\author{Andr\'{e} Mas\thanks{Institut de Mod\'{e}lisation Math\'{e}matique de
Montpellier, CC 051, Universit\'{e} Montpelllier 2, place Eug\`{e}ne
Bataillon, 34095 Montpellier Cedex 5, France, mas@math.univ-montp2.fr}\\Universit\'{e} Montpellier 2}
\title{Local Functional Principal Component Analysis\bigskip}
\date{}
\maketitle

\begin{abstract}
Covariance operators of random functions are crucial tools to study the way
random elements concentrate over their support. The principal component
analysis of a random function $X$ is well-known from a theoretical viewpoint
and extensively used in practical situations. In this work we focus on local
covariance operators. They provide some pieces of information about the
distribution of $X$ around a fixed point of the space $x_{0}$. A description
of the asymptotic behaviour of the theoretical and empirical counterparts is
carried out. Asymptotic developments are given under assumptions on the
location of $x_{0}$ and on the distributions of projections of the data on the
eigenspaces of the (non-local) covariance operator.

\end{abstract}

\section{Introduction}

\subsection{The general framework}

A recent and considerable interest has been given along the last years to the
statistical analysis for functional data. Usually mathematical statistics or
probability theory deal with data modelled as random variables i.e. measurable
mappings from an abstract probability space $\left(  \Omega,\mathcal{A}%
,\mathbb{P}\right)  $ to a finite dimensional space. If a sample is denoted
$X_{1},..,X_{n},$ the $X$'s take values classically in $\mathbb{R}$ or
$\mathbb{R}^{p}$. In our framework we turn to random elements with values in
an infinite-dimensional function space denoted $\mathcal{F}$. In the sequel
$\mathcal{F}$ will be endowed at least with a separable Banach space structure
with norm $\left\Vert \cdot\right\Vert $. In everyday's life, situation where
such data appear are quite common : monitoring the value of a share yields a
random curve $X\left(  t\right)  $ where $t$ runs along the quotation time.
Even more basically, observing the temperature at a given place along the day
and during $n$ consecutive days provides a "theoretical" sample $X_{1}\left(
t\right)  ,...,X_{n}\left(  t\right)  $ where $t\in\left[  0,24\right]  .$
Here "theoretical" means that temperatures will be recorded each day at fixed
moments. For instance $X_{1}$ will be observed at $t_{1},...t_{m_{1}},$ hence
the true curve $X_{1}\left(  t\right)  $ will have to be reconstructed from
the finite real valued sample $X_{1}\left(  t_{1}\right)  ,...,X_{1}\left(
t_{m_{1}}\right)  $ by interpolation techniques such as splines, wavelets,
cosine bases, etc. For further references about this topic see Chui (1992) and
Antoniadis, Oppenheim (1995) about wavelets and de Boor (1978) ou Dierckx
(1993), about splines.

It turns out that probabilists have studied such random elements for a much
longer time than statisticians (first works on the Brownian motion date back
to the XIX$^{th}$ century), the first monograph dedicated to functional data
was published in 1991 (Ramsay and Silverman (1991)). Modern computers make it
now possible to carry out calculations for very high dimensional vectors and
practical statisticians have shifted their interests to functional data or to
so-called "high-dimensional problems" that fall within the scope of this
paper. However a theoretical gap remains because not all asymptotic results
have been given yet for such data. Besides probability theory unfortunately
sometimes just give clues and not solutions to typically statistical issues
(see later the paragraph devoted to small ball problems).

The interested reader could get familiar with the applied aspects by reading
the monographs by Ramsay and Silverman or by Ferraty and Vieu (2006). Many
probabilistic results will be found in Vakhania, Tarieladze and Chobanian
(1987) or in Ledoux and Talagrand (1991).

In this setting we can define the expectation of $X$ as :%
\[
\mathbb{E}X=\int_{\Omega}Xd\mathbb{P=}\int_{\mathcal{F}}xd\mu_{X}\left(
x\right)  \in\mathcal{F}%
\]
where $\mu_{X}=\mathbb{P\circ}X^{-1}$ is the image of measure of $\mathbb{P}$
through the mapping $X$. The integral is of Bochner type and is defined
whenever the real valued random variable $\left\Vert X\right\Vert $ is integrable.

In this article we focus on a very useful and common statistical technique :
principal component analysis (PCA for short). The functional version of the
PCA was initially studied by Dauxois, Pousse and Romain (1982). We refer to
this seminal article for a complete mathematical definition. Briefly speaking
functional PCA of a process $X\left(  \cdot\right)  $ comes down to the
spectral analysis of the covariance operator associated to the sample (see
below for definitions). We refer for instance to Silverman (1996), Oca\~{n}a,
Aguilera and Valderrama (1999), Kneip and Utikal (2001), Yao, Muller and Wang
(2005), or Cardot, Mas and Sarda (2007) to overview some applications and
extensions of the functional PCA. Also note that He, Muller and Wang (2003)
introduced a version of canonical analysis for random functions.

\subsection{The Hilbert setting}

The abstract framework defined above is too general for statistical purposes.
We will restrict ourselves to special spaces $\mathcal{F}$ but we will gain in
terms of interpretation of our assumptions and results. Indeed we will assume
once and for all that $\mathcal{F}=H$ is a separable Hilbert space endowed
with inner product $\left\langle \cdot,\cdot\right\rangle $ and associated
norm $\left\Vert \cdot\right\Vert .$ Several reasons can explain this
tightening. First of all we know that each function $X$ may be decomposed in
denumerable bases. In practical situations these bases enable to get
observation-curves from the discretized ones (see first paragraph above) by
interpolation methods mentioned above. The following step is then achieved :%
\[
\left[  X_{1}\left(  t_{1}\right)  ,...,X_{1}\left(  t_{m_{1}}\right)
\right]  \Rightarrow X_{1}\left(  t\right)  =\sum_{k=1}^{N}c_{k}e_{k}\left(
t\right)  .
\]

The Sobolev spaces $W^{m,2}$ are classical examples of such Hilbert spaces :%
\[
W^{m,2}\left(  \left[  0,T\right]  \right)  =\left\{  f\in L^{2}\left(
\left[  0,T\right]  \right)  :\sum_{k=0}^{m}\int_{0}^{T}\left(  f^{\left(
k\right)  }\left(  s\right)  \right)  ^{2}ds<+\infty\right\}
\]
where $T>0$ and $f^{\left(  k\right)  }$ denotes here the derivative of order
$k$ of $f$.

Besides as we will focus on covariance operators the Hilbert setting yields
considerable simplifications. Bounded linear operators were extensively
studied as well as their spectral properties (see references below). We
introduce the two following operator spaces and associated norms.

The Banach space $\mathcal{L}\left(  H,H\right)  =\mathcal{L}$ is the
classical space of bounded operator endowed with the norm defined for each $T$
in $\mathcal{L}$ by :%
\[
\left\Vert T\right\Vert _{\infty}=\sup_{x\in\mathcal{B}_{1}}\left\Vert
Tx\right\Vert ,
\]
where $\mathcal{B}_{1}$ is the unit sphere of $H.$ The Hilbert space
$\mathcal{L}_{2}$ is the space of Hilbert-Schmidt operators, ($\mathcal{L}%
_{2}\subset\mathcal{L}$) i.e. the space of those operators $T$ such that,
given a basis of $H,$ say $\left(  e_{k}\right)  _{k\in\mathbb{N}},$%
\[
\left\Vert T\right\Vert _{2}=\sum_{k=1}^{+\infty}\left\Vert Te_{k}\right\Vert
^{2}<+\infty.
\]
It is a well-known fact that $\mathcal{L}_{2}$ is a separable Hilbert space
whenever $H$ is. The inner product in $\mathcal{L}_{2}$ is :%
\[
\left\langle T,S\right\rangle _{2}=\sum_{k=1}^{+\infty}\left\langle
Te_{k},Se_{k}\right\rangle
\]
and does not depend on the choice of the basis $\left(  e_{k}\right)
_{k\in\mathbb{N}}$. The space of trace-class (or nuclear) operators
$\mathcal{L}_{1}$ endowed with norm $\left\Vert \cdot\right\Vert _{1}$ will be
mentioned sometimes in the paper. These norms are not equivalent and%
\[
\left\Vert \cdot\right\Vert _{\infty}\leq\left\Vert \cdot\right\Vert _{2}%
\leq\left\Vert \cdot\right\Vert _{1}.
\]
The canonical injections from $\mathcal{L}_{1}$ onto $\mathcal{L}_{2}$ and
from $\mathcal{L}_{2}$ onto $\mathcal{L}$ are consequently continuous. For
further information on linear operators we refer to Schmeidler (1965), Weidman
(1980), Dunford-Schwartz (1988), Gohberg, Goldberg and Kaashoek (1991) amongst
many others.

\section{Covariance and local covariance operators of random Hilbert elements}

Let the tensor product between $u$ and $v$ in $H$ stand for the one-rank
operator from $H$ to $H$ by :%
\[
\left(  u\otimes v\right)  \left(  t\right)  =\left\langle u,t\right\rangle v
\]
for all $t$ in $H.$

Since covariance operators are undern concern within functional PCA, we should
first of all define them and give some of their main features. The theoretical
covariance operator $\Gamma$ and its empirical counterpart, $\Gamma_{n}$,
based on the independent and identically distributed sample $X_{1},...,X_{n}$
are symmetric positive trace class operators from $H$ to $H$ defined by :%
\begin{align}
\Gamma &  =\mathbb{E}\left(  \left(  X_{1}-\mathbb{E}X_{1}\right)
\otimes\left(  X_{1}-\mathbb{E}X_{1}\right)  \right)  ,\label{thcovop}\\
\Gamma_{n}  &  =\dfrac{1}{n}\sum_{k=1}^{n}\left(  X_{k}-\overline{X}%
_{n}\right)  \otimes\left(  X_{k}-\overline{X}_{n}\right)  \label{empcovop}%
\end{align}
where%
\[
\overline{X}_{n}=\frac{1}{n}\sum_{k=1}^{n}X_{k}.
\]
If $X$ is centered $\mathbb{E}X=0$ and $\Gamma=\mathbb{E}\left(  X\otimes
X\right)  .$ By $\left(  \lambda_{k},e_{k}\right)  $ we denote the $k^{th}$
eigenelements (eigenvalues/eigenvectors) of $\Gamma$. The $\lambda_{k}$'s are
positive and we set $\lambda_{1}\leq\lambda_{2}\leq...,$ and $\left(
\lambda_{k}\right)  _{k\in\mathbb{N}}\in l_{1}$.

In the Hilbert setting, the distribution of a centered random element say $X$
may be characterized a very simple way. Indeed if $=_{d}$ denotes equality in
distribution :%
\begin{equation}
X=_{d}\sum_{k=1}^{+\infty}\sqrt{\lambda_{k}}\xi_{k}e_{k} \label{KL}%
\end{equation}
where the $\xi_{k}$'s are centered non-correlated real random variables with
unit variance. The above decomposition is often referred to as the
Karhunen-Lo\`{e}ve development or development of $X$ with respect to its
reproducing kernel Hilbert space (RKHS). For definition and studies of RKHS we
refer to Berlinet, Thomas-Agnan (2004).

Now we introduce the definition of local covariance operator. Let $K$ be a
kernel, that is a nonnegative function defined on $\mathbb{R}^{+}$, and such
that $\int K\left(  s\right)  ds=1$ and let $h=h\left(  n\right)  $ be a
bandwidth. Pick $x_{0}$ a fixed vector in $H$ ($x_{0}$ will be a function if
$H$ is a space of functions). In the sequel we will often need to decompose
$x_{0}$ in the basis $e_{k}$ :%
\[
x_{0}=\sum_{k=1}^{+\infty}\left\langle x_{0},e_{k}\right\rangle e_{k}%
\]
where%
\[
\sum_{k=1}^{+\infty}\left\langle x_{0},e_{k}\right\rangle ^{2}<+\infty.
\]

\begin{definition}
The theoretical local covariance operator of $X$ at $x_{0}\in H$ based on the
kernel $K$ and its empirical counterpart are respectively defined by :%
\begin{align}
\Gamma_{K}  &  =\mathbb{E}\left(  K\left(  \frac{\left\Vert X_{1}%
-x_{0}\right\Vert }{h}\right)  \left(  \left(  X_{1}-x_{0}\right)
\otimes\left(  X_{1}-x_{0}\right)  \right)  \right)  ,\label{lcopth}\\
\Gamma_{K,n}  &  =\dfrac{1}{n}\sum_{k=1}^{n}K\left(  \frac{\left\Vert
X_{k}-x_{0}\right\Vert }{h}\right)  \left(  \left(  X_{k}-x_{0}\right)
\otimes\left(  X_{k}-x_{0}\right)  \right)  . \label{lcopemp}%
\end{align}

\end{definition}

Even if $\Gamma_{K}$ and $\Gamma_{K,n}$ may have already been introduced
elsewehere in articles dealing with functional data (see for instance Ferraty,
Mas, Vieu (2007)), it is the first attempt, up to the author's knowledge, to
provide these operators with a name. Note that $\Gamma_{K}$ implicitely
depends on $n$ through $h,$ even if this index does not explicitely appears.
These operators are crucial in the nonparametric estimation of the regression
function by local linear methods amongst others (see the conclusion for
further details).

\begin{proposition}
The local covariance operators $\Gamma_{K}$ and $\Gamma_{K,n}$ are positive,
selfadjoint. Besides they are trace class whenever $K$ is bounded and%
\[
\mathbb{E}\left(  \left\Vert X_{1}\right\Vert ^{2}\right)  <+\infty
\]
or when $K$ is bounded and compactly supported.
\end{proposition}

The proposition is plain since the trace-class norm of $\Gamma_{K}$ for
instance is bounded by :
\begin{align*}
&  \mathbb{E}\left\Vert K\left(  \frac{\left\Vert X_{1}-x_{0}\right\Vert }%
{h}\right)  \left(  \left(  X_{1}-x_{0}\right)  \otimes\left(  X_{1}%
-x_{0}\right)  \right)  \right\Vert _{1}\\
&  =\mathbb{E}\left[  K\left(  \frac{\left\Vert X_{1}-x_{0}\right\Vert }%
{h}\right)  \left\Vert X_{1}-x_{0}\right\Vert ^{2}\right]  .
\end{align*}
In the sequel by $\left(  \lambda_{k}\right)  _{k\in\mathbb{N}}$ (resp.
$\left(  \lambda_{k,n}\right)  _{k\in\mathbb{N}}$) and $\left(  \pi
_{k}\right)  _{k\in\mathbb{N}}$ (resp. $\left(  \pi_{k,n}\right)
_{k\in\mathbb{N}}$) we denote the eigenvalues and the associated eigenvectors
of $\Gamma_{K}$ (resp. $\Gamma_{K,n}$)

The main goal of this paper is to describe the asymptotic behaviour of
$\Gamma_{K}$ and $\Gamma_{K,n}$ and to derive results for their eigenelements.

\section{Intermezzo about Gamma varying functions and the small ball problem}

It may be proved that whenever the random variable $X$ takes values in a
finite dimensional space, say $\mathbb{R}^{p}$, the covariance operators
defined at displays (\ref{lcopth}) and (\ref{lcopemp}) depend on the value of
the density of $X,$ if we assume that $X$ admits a density at point $x_{0}$.
Besides one may use a larger class of kernels (here since $\left\Vert
X_{1}-x_{0}\right\Vert $ is positive we are restricted to kernels with
positive support which damages the rates of convergence). It is simple to
prove that when $X$ is a real-valued random variable that exhibits a non-null
density $f$ at $x_{0}$ with some regularity around $x_{0}$ :%
\[
\mathbb{E}\left[  K\left(  \frac{X_{1}-x_{0}}{h}\right)  \left(  X_{1}%
-x_{0}\right)  ^{2}\right]  \sim h^{3}f_{X}\left(  x_{0}\right)  \int
u^{2}K\left(  u\right)  du.
\]
We refer to Fan (1993) for illustrating the issues of asymptotics for
truncated moments in nonparametric regression estimation.

In our infinite dimensional framework the situation is quite different : since
Lebesgue's measure is not defined on Hilbert spaces, the notion of 'density'
cannot be defined either. It turns out that the density will be replaced by
the "small ball probability of $X$", which is nothing but the cumulative
density function for the norm of $X$ (or $X-x_{0}$) in a neighborhood of $0$
i.e. $\mathbb{P}\left(  \left\Vert X-x_{0}\right\Vert <\varepsilon\right)  $
for $\varepsilon\downarrow0$. The random variable $X$ may be replaced by any
process $Z_{t},t\in T$. The study of small ball probabilites is not new and is
connected with the theory of large deviations (since $\left\Vert
X-x_{0}\right\Vert /\varepsilon$ will be large when $\varepsilon$ decays to
zero). We refer to Ledoux-Talagrand (1991), Li, Linde (1999), Li, Shao (2001)
for some more information about this topic.

Since these small ball probabilities appear in the main results of this
article and are of much importance within the proofs, we collect a few results
about them in order to be more illustrative. We provide two examples directly
based on display (\ref{KL}). The small probability then heavily depends on the
rate of decay of the eigenvalues of the covariance operator, namely the
$\lambda_{k}$'s.

If the rate of decay is arithmetic $\lambda_{k}\asymp k^{-\left(
1+\alpha\right)  }$. The problem was solved in Mayer-Wolf, Zeitouni (1993).
They get :%
\begin{equation}
\mathbb{P}\left(  \left\Vert X\right\Vert <\varepsilon\right)  \sim\exp\left(
-\frac{C\left(  \alpha\right)  }{\varepsilon^{1/\alpha}}\right)  . \label{PB1}%
\end{equation}
where $C\left(  \alpha\right)  $ is some positive constant.

When the rate of decay is exponential : $\lambda_{k}\asymp\exp\left(
-ak\right)  ,$ the calculations may not have been carried out. I did not find
them elsewhere. They may be derived from formula (10) in Dembo, Mayer-Wolf,
Zeitouni (1995).

\begin{proposition}
\label{Juliet}When $\lambda_{k}\asymp\exp\left(  -ak\right)  $ in (\ref{KL}),
then for $\varepsilon\rightarrow0$ :%
\begin{equation}
\mathbb{P}\left(  \left\Vert X\right\Vert <\varepsilon\right)  \sim\sqrt
{\frac{\alpha}{-\pi\log\left(  \varepsilon\right)  }}\exp\left(  -\frac
{1}{4\alpha}\left[  \log\left(  \varepsilon\right)  \right]  ^{2}\right)
\label{PB2}%
\end{equation}

\end{proposition}

The proof of this apparently new result is postponed to the end of the last section.

Let us leave the small ball probability for a moment. At this point we need to
give some properties of a class of real functions. The statistician may be
familiar with the definition of functions with regular variations since they
appear in the theory of extremes. For instance $f:\mathbb{R}\rightarrow
\mathbb{R}$ is regularly varying at $0$ with index $d$ if , for all fixed $x$
in $\mathbb{R}$ :%
\[
\lim_{h\rightarrow0}\frac{f\left(  hx\right)  }{f\left(  h\right)  }=x^{d}.
\]
A less known class of functions studied in the theory of regular variations is
the so-called "class $\Gamma$". This class $\Gamma$ will be of much use in the
sequel. It was introduced by de Haan (1971), see also de Haan (1974) in
connection with the theory of extemes. But Ga\"{\i}ffas (2005) used it to
model the distribution of "rare" inputs in a non-parametric regression model :
a density which is null and $\Gamma$-varying at $x_{0}$ will generate a
distribution which rarely visit $x_{0}$.

\begin{definition}
A function $f$ belongs to the class $\Gamma$ at $0$ (or is $\Gamma$-varying at
$0$) if there exists a measurable positive function $\rho$ such that for all
$x\in\mathbb{R}$ :%
\begin{equation}
\lim_{h\rightarrow0^{+}}\frac{f\left(  h+\rho\left(  h\right)  x\right)
}{f\left(  h\right)  }=\exp\left(  x\right)  . \label{defgamma}%
\end{equation}
The function $\rho$ is called the auxiliary function of $f$.
\end{definition}

We refer to Chapter 3.10 in Bingham, Goldie, Teugels (1987) for a deeper
presentation and the essential properties of regularly varying functions and
of the class $\Gamma$ (see p.174-180). The reader should be aware that, in
this book, the authors consider only functions that are $\Gamma$-varying at
infinity. Their definitions and properties must be adapted to our setting :
here functions are $\Gamma$-varying at $0$. We collect now only those
properties which will be used throughout the proofs :

\textbf{Fact 1} : If $f\in\Gamma,$ for all $x\in\left[  0,1\right[  $%
\begin{equation}
\lim_{h\rightarrow0^{+}}\frac{f\left(  hx\right)  }{f\left(  h\right)  }=0.
\label{F1}%
\end{equation}

\textbf{Fact 2} : If $\rho$ is the auxiliary function of $f\in\Gamma$, then%
\begin{align}
&  \frac{\rho\left(  s\right)  }{s}\underset{s\rightarrow0}{\rightarrow
}0,\label{F2}\\
&  \frac{\rho\left(  s+x\rho\left(  s\right)  \right)  }{\rho\left(  s\right)
}\underset{s\rightarrow0}{\rightarrow}1 \label{F2bis}%
\end{align}
when $s$ goes to $0$ and for all $x\in\mathbb{R}$.

\textbf{Fact 3} : If $f\in\Gamma$ then $F\left(  x\right)  =\int_{0}%
^{x}f\left(  s\right)  ds$ belongs to the class $\Gamma$ too and%
\begin{equation}
\int_{0}^{h}f\left(  s\right)  ds\underset{h\rightarrow0}{\sim}f\left(
h\right)  \rho\left(  h\right)  . \label{F3}%
\end{equation}

Now we turn again to the small ball probabilities. The following Proposition
explains why the class $\Gamma$ was introduced.

\begin{proposition}
\label{MED}Functions defined on displays (\ref{PB1}) and (\ref{PB2}) are both
$\Gamma$-varying at $0$ with auxiliary functions~:%
\[
\rho\left(  s\right)  =\frac{\alpha}{C\left(  \alpha\right)  }s^{1+1/\alpha}%
\]
and
\[
\rho\left(  s\right)  =-s/\left(  2\alpha\log s\right)
\]
respectively.
\end{proposition}

The proof is omitted since it is straightforward. The previous Proposition is
quite important for the sequel.

It is seen from (\ref{lcopth}) for instance that the random element $X$ is
shifted from the origin by $-x_{0}.$ Obviously small ball probabilites defined
at displays (\ref{PB1}) or (\ref{PB2}) do not exactly match our goals. The
shift, $-x_{0}$ is nonrandom but however the small balls probabilities may
tremendously differ from those given above. It turns out that when $x_{0}$
lies in the reproducing kernel Hilbert space of $X$ :%
\begin{equation}
\mathbb{P}\left(  \left\Vert X-x_{0}\right\Vert <\varepsilon\right)  \sim
C\left(  x_{0}\right)  \mathbb{P}\left(  \left\Vert X\right\Vert
<\varepsilon\right)  \label{craig}%
\end{equation}
where $C\left(  x_{0}\right)  $ is a constant which depends only on $x_{0}$
and $\asymp$ should replace $\sim$ in (\ref{PB1}) and (\ref{PB2}). The
articles by Li, Linde (1993) or Kuelbs, Li, Linde (1994) deal with the small
ball problem for shifted balls.Since the small ball probability of $X$ appears
explicitely in the main results of this paper it is denoted for simplicity :%
\[
F_{x_{0}}\left(  \varepsilon\right)  =F\left(  \varepsilon\right)
=\mathbb{P}\left(  \left\Vert X-x_{0}\right\Vert <\varepsilon\right)  .
\]
Within the proofs several calculations must be carried out that involve some
analytic properties of $F$ and we need to announce the following claim,
inspired by Proposition \ref{MED} and (\ref{craig})

\begin{claim}
The small ball probability functions $F$ of shifted random elements in $H$
(here $X-x_{0}$) are assumed to belong to the class $\Gamma.$
\end{claim}

The assumptions needed to define correctly $F_{x_{0}}\left(  \varepsilon
\right)  $ that may appear in Li, Linde (1993) or Kuelbs, Li, Linde (1994) are
supposed to hold in addition to those that will be given below.

\begin{remark}
This claim directly leads us to arising the following question : is the small
probability (at $0$) of any process defined by (\ref{KL}) $\Gamma$-varying at
zero ? Answering yes would provide a universal "representation" of these small
ball probability functions (see Theorem 3.10.8 p.178 in Bingham, Goldie,
Teugels (1987)). At this point we can answer only in the two important special
cases mentioned above but this issue is under investigation (see Mas (2007)).
\end{remark}

\section{Main results}

We study convergence for random (or not) operators. This section is tiled into
three subsections. In the first one we provide theorems dealing with
asymptotics for the cells of some infinite matrix. Exact constants are
computed. In the second subsection we get bounds in supremum or
Hilbert-Schmidt norm for the operator(s) under concern. The third deals with
the empirical local covariance operator $\Gamma_{K,n}.$

\subsection{Cell-by-cell results}

First let us introduce the assumptions needed in the sequel.\medskip

\textbf{Assumption }$\mathbf{A}_{1}$\textbf{ :} \textit{There exists a basis,
say }$e_{p}$\textit{ in which the finite dimensional distributions of }%
$X,$\textit{ the }$\left\langle X,e_{p}\right\rangle $\textit{'s are
independent and all have a density }$f_{p}$\textit{. This density is such that
}$\left(  f_{p}\right)  ^{\left(  i\right)  }\left(  \left\langle x_{0}%
,e_{p}\right\rangle \right)  \neq0$\textit{ for }$i\in\left\{  0,1,2\right\}
.$ \textit{Besides the density of the nonnegative real variable}%
\[
\sqrt{\sum_{k=1}^{+\infty}\left\langle X-x_{0},e_{k}\right\rangle ^{2}%
}=\left\Vert X-x_{0}\right\Vert
\]
\textit{exists in a neighborhood of }$0$ and \textit{belongs to the class
}$\Gamma$\textit{ with auxiliary function }$\rho$\textit{.}

\begin{remark}
The notation $e_{p}$ should not be misleading. The basis involved in
Assumption $\mathbf{A}_{1}$ needs not to be the basis of eigenvalues of the
operator $\Gamma.$ But since in the important case of a gaussian random
element $X$ -with eigenvalues decaying at an arithmetic or geometric rate-
$\mathbf{A}_{1}$ always holds for this special basis we will abusively keep
the same notation.\medskip
\end{remark}

Assuming that the finite dimensional distributions are independent is needed
to alleviate the proofs and to get exact constant in asymptotic expansions.
Milder hypotheses on the joint distribution of the couple $\left(
\left\langle X,e_{k}\right\rangle ,\left\Vert X\right\Vert \right)  $ could
certainly prevail at the expense of more tedious caclulations as will be seen
from the proofs.

\begin{remark}
Within $\mathbf{A}_{1}$ the assumption $\left(  f_{p}\right)  ^{\left(
i\right)  }\left(  \left\langle x_{0},e_{p}\right\rangle \right)  \neq
0$\textit{ for }$i\in\left\{  0,1,2\right\}  $ could be replaced by the more
general one : "Let us denote%
\begin{align*}
N_{p}^{0}  &  =\inf\left\{  k:f_{p}^{\left(  2k\right)  }\left(  \left\langle
x_{0},e_{p}\right\rangle \right)  \neq0\right\}  ,\\
N_{p}^{1}  &  =\inf\left\{  k:f_{p}^{\left(  2k+1\right)  }\left(
\left\langle x_{0},e_{p}\right\rangle \right)  \neq0\right\}  ,\\
N_{p}^{2}  &  =\inf\left\{  k>N_{p}^{0}:f_{p}^{\left(  2k\right)  }\left(
\left\langle x_{0},e_{p}\right\rangle \right)  \neq0\right\}
\end{align*}
and assume that $N_{p}^{1}$ and $N_{p}^{2}$ are finite for all $p$". But once
more we prefer to lose generality and gain readability. Also note that
switching Assumption 1 to the one involving the $N_{p}^{k}$'s leads to
modified results in Theorem \ref{LCO} : indeed the speed of convergence would
then depend on $N_{p}^{0},$ $N_{p}^{1}$ and $N_{p}^{2}$.\medskip
\end{remark}

\textbf{Assumption }$\mathbf{A}_{\mathbf{2}}$\textbf{ :} \textit{The kernel
}$K$\textit{ is bounded, }$\left[  0,1\right]  $\textit{-supported, }$K\left(
1\right)  >0$ and
\[
\sup_{s\in\left[  0,1\right]  }\left\vert K^{\prime}\left(  s\right)
\right\vert <+\infty
\]

This assumption is not too restrictive and could certainly be replaced by a
milder one. But it is out of the scope of this article to provide minimal
conditions on the kernel $K$.\medskip

We start with a development of $\Gamma_{K}$. Let $\delta_{i,j}$ be the
Kronecker symbol ($\delta_{i,j}=1$ if and only if $i=j,$ $0$ otherwise).

\begin{theorem}
\label{LCO}Assume $\mathbf{A}_{1}$ and $\mathbf{A}_{\mathbf{2}}$. When $h$
goes to $0$ the operator $\Gamma_{K}$ also tends to zero. And the following
holds : for fixed $i$ and $j$ in $\mathbb{N}$,%
\begin{equation}
\left\langle \Gamma_{K}\left(  h\right)  \left(  e_{i}\right)  ,e_{j}%
\right\rangle \sim v\left(  h\right)  \delta_{i,j}+w\left(  h\right)
\mathcal{R}_{ij}, \label{Clint}%
\end{equation}
where $v\left(  h\right)  $ and $w\left(  h\right)  $ are two real nonnegative
sequences defined by :%
\begin{align}
v\left(  h\right)   &  =\mathbb{E}\left(  K\left(  \frac{\left\Vert
X-x_{0}\right\Vert }{h}\right)  \left\Vert X-x_{0}\right\Vert \rho\left(
\left\Vert X-x_{0}\right\Vert \right)  \right)  ,\label{V}\\
w\left(  h\right)   &  =\mathbb{E}\left(  K\left(  \frac{\left\Vert
X-x_{0}\right\Vert }{h}\right)  \left\Vert X-x_{0}\right\Vert ^{2}\rho
^{2}\left(  \left\Vert X-x_{0}\right\Vert \right)  \right)  \label{W}%
\end{align}
and where the doubly indexed field $\mathcal{R}$ is defined by :%
\begin{align*}
\mathcal{R}_{ii}  &  =\frac{f_{i}^{\prime\prime}\left(  \left\langle
x_{0},e_{i}\right\rangle \right)  }{f_{i}\left(  \left\langle x_{0}%
,e_{i}\right\rangle \right)  },\\
\mathcal{R}_{ij}  &  =\frac{f_{i}^{\prime}\left(  \left\langle x_{0}%
,e_{i}\right\rangle \right)  }{f_{i}\left(  \left\langle x_{0},e_{i}%
\right\rangle \right)  }\frac{f_{j}^{\prime}\left(  \left\langle x_{0}%
,e_{j}\right\rangle \right)  }{f_{j}\left(  \left\langle x_{0},e_{j}%
\right\rangle \right)  },\quad i\neq j.
\end{align*}

\end{theorem}

Theorem \ref{LCO} provides asymptotics for each cell of the infinite
dimensional matrix $\Gamma_{K}$ when expressed in the basis $\left(
e_{i}\right)  _{1\leq i\leq n}.$ Introducing the operator $\mathcal{R}$
defined in the basis $\left(  e_{i}\right)  _{1\leq i\leq n}$ by $\left\langle
\mathcal{R}e_{i},e_{j}\right\rangle =\mathcal{R}_{ij}$ we could rephrase this
theorem by saying that "$\Gamma_{K}$ is asymptotically equivalent "cell by
cell" to the operator $v\left(  h\right)  I+w\left(  h\right)  \mathcal{R}$".

At this point the reader is not given much information on both sequences
$v\left(  h\right)  $ and $w\left(  h\right)  .$ It is actually basic to see
that both sequences tend to zero. The next subsection will provide the reader
with a more explicit description of the rate of decrease.

The next Proposition and the two next remarks give seminal properties of
operator $\mathcal{R}$.

\begin{proposition}
\label{nova}If both following conditions hold%
\begin{align}
\mathbf{C}_{1}  &  :\sum_{i=1}^{+\infty}\left(  \frac{f_{i}^{\prime\prime
}\left(  \left\langle x_{0},e_{i}\right\rangle \right)  }{f_{i}\left(
\left\langle x_{0},e_{i}\right\rangle \right)  }\right)  ^{2}<+\infty
,\label{g1}\\
\mathbf{C}_{2}  &  :\sum_{i=1}^{+\infty}\left(  \frac{f_{i}^{\prime}\left(
\left\langle x_{0},e_{i}\right\rangle \right)  }{f_{i}\left(  \left\langle
x_{0},e_{i}\right\rangle \right)  }\right)  ^{2}<+\infty, \label{g1bis}%
\end{align}
$\mathcal{R}$ is Hilbert-Schmidt. If (\ref{g1}) is replaced with
\begin{equation}
\mathbf{C}_{1}^{\prime}:\left(  \frac{f_{i}^{\prime\prime}\left(  \left\langle
x_{0},e_{i}\right\rangle \right)  }{f_{i}\left(  \left\langle x_{0}%
,e_{i}\right\rangle \right)  }\right)  _{i\in\mathbb{N}}\in l_{\infty}%
\quad\left(  \mathrm{resp\ }c_{0}\right)  \label{g2}%
\end{equation}
the operator $\mathcal{R}$ is bounded (resp. compact).\newline When either
(\ref{g1bis}) or (\ref{g2}) does not hold, $\mathcal{R}$ is a symmetric
unbounded operator.
\end{proposition}

The proof of this Proposition is omitted since it is a consequence of the
following remark.

\begin{remark}
The operator $\mathcal{R}$ may be rewritten :%
\[
\mathcal{R}=\tau\otimes\tau+\mathrm{diag}\left(  s_{i}\right)
\]
where%
\begin{align*}
\tau &  =\left(  \frac{f_{1}^{\prime}\left(  \left\langle x_{0},e_{1}%
\right\rangle \right)  }{f_{1}\left(  \left\langle x_{0},e_{1}\right\rangle
\right)  },\frac{f_{2}^{\prime}\left(  \left\langle x_{0},e_{2}\right\rangle
\right)  }{f_{2}\left(  \left\langle x_{0},e_{2}\right\rangle \right)
},...\right) \\
&  =\left(  \left(  \ln f_{1}\right)  ^{\prime}\left(  \left\langle
x_{0},e_{1}\right\rangle \right)  ,\left(  \ln f_{2}\right)  ^{\prime}\left(
\left\langle x_{0},e_{2}\right\rangle \right)  ,...\right)  ,\\
s_{i}  &  =\left(  \ln f_{i}\right)  ^{\prime\prime}\left(  \left\langle
x_{0},e_{i}\right\rangle \right)  ,
\end{align*}
and $\mathrm{diag}\left(  s_{i}\right)  $ denotes a diagonal operator
expressed in the basis $e_{i}$ with $i^{th}$ term $s_{i}$. When (\ref{g1bis})
holds $\tau\in H$.
\end{remark}

Before going into deeper details we should examine a typical situation, namely
the case when $X$ is a gaussian random element. The following Proposition
shows that even in this basic situation, serious problems occur.

\begin{proposition}
If $X$ is gaussian and centered :%
\begin{align*}
\frac{f_{i}^{\prime}\left(  \left\langle x_{0},e_{i}\right\rangle \right)
}{f_{i}\left(  \left\langle x_{0},e_{i}\right\rangle \right)  }  &
=-\frac{\left\langle x_{0},e_{i}\right\rangle }{\lambda_{i}},\\
\frac{f_{i}^{\prime\prime}\left(  \left\langle x_{0},e_{i}\right\rangle
\right)  }{f_{i}\left(  \left\langle x_{0},e_{i}\right\rangle \right)  }  &
=\left(  \frac{\left\langle x_{0},e_{i}\right\rangle }{\lambda_{i}}\right)
^{2}-\frac{1}{\lambda_{i}}%
\end{align*}
and conditions (\ref{g1bis}) and (\ref{g2}) of Proposition \ref{nova} cannot
hold together which also means that the operator $\mathcal{R}$ is always
unbounded in this setting.
\end{proposition}

Let us focus on these conditions since they may be easily interpreted. Indeed
assuming (\ref{g1bis}) :%
\[
\sum_{i=1}^{+\infty}\left(  \frac{\left\langle x_{0},e_{i}\right\rangle
}{\lambda_{i}}\right)  ^{2}<+\infty
\]
means that the coordinates of $x_{0}$ should decrease much quicker to zero
than the eigenvalues which also means that $x_{0}$ should be smoother (more
regular) than $X$ itself. On the other hand (\ref{g2}) will hold whenever
$\left\langle x_{0},e_{i}\right\rangle ^{2}=\lambda_{i}+\lambda_{i}^{2}%
\tau_{i}$ (where $\tau\in l_{\infty}\quad\left(  \mathrm{resp\ }c_{0}\right)
$), and $x_{0}$ should then be as regular as $X$ but -surprisingly- not more...

Considering again the unusual conditions $\mathbf{C}_{1},\mathbf{C}_{2}$ the
reader may become suspicious and the situation in the gaussian framework makes
it legitimate to wonder whether there exists a family of densities and an
$x_{0}$ such that these conditions hold. The answer is positive and gives
birth to the following Proposition.

\begin{proposition}
Let $f_{i}$ be the symmetric density defined on $\mathbb{R}$ by :%
\[
f_{i}\left(  x\right)  =\frac{6}{3^{6}}\frac{1}{\lambda_{i}^{2}}\left(
27\lambda_{i}^{3/2}-\left\vert x\right\vert ^{3}\right)  1\mathrm{I}_{\left\{
x\leq3\lambda_{i}^{1/2}\right\}  }%
\]
and take $\left\langle x_{0},e_{i}\right\rangle =x_{i}$ such that
\[
\sum_{i=1}^{+\infty}\frac{x_{i}^{2}}{\lambda_{i}^{3}}<+\infty.
\]
Then (\ref{g1}) and (\ref{g1bis}) both hold.
\end{proposition}

The proof of the Proposition is omitted since it stems from straightforward
computations.\medskip

\subsection{Norm results for the local covariance operator}

The next issue is obviously to strengthen Theorem \ref{LCO} : Is it possible
to replace the "cell by cell" or "componentwise" convergence by convergence in
norm ? First of all note that we may expect the rate of convergence to be
$v\left(  h\right)  $ but we have to be cautious for several reasons :

\begin{itemize}
\item First of all $\mathcal{R}$ may be unbounded. In that situation we cannot
expect a result such as :%
\[
\Gamma_{K}-\left\{  v\left(  h\right)  I+w\left(  h\right)  \mathcal{R}%
\right\}  \rightarrow0
\]
in norm since $w\left(  h\right)  \mathcal{R}$ may not even be bounded whereas
$\Gamma_{K}$ is.

\item We may have $\Gamma_{K}-v\left(  h\right)  I\rightarrow0$ but we cannot
get%
\[
\left\Vert \frac{\Gamma_{K}}{v\left(  h\right)  }-I\right\Vert _{\infty
}\rightarrow0
\]
for topological reasons : $\frac{\Gamma_{K}}{v\left(  h\right)  }$ is for all
$h$ a compact operator and cannot converge to the identity operator (which is
not compact) since $\mathcal{L}_{c}$ is a closed subspace of $\mathcal{L}$.

\item Even worse : $\frac{\Gamma_{K}}{v\left(  h\right)  }$ may be
asymptotically bounded or unbounded i.e.
\[
\limsup_{h\rightarrow0}\left\Vert \frac{\Gamma_{K}}{v\left(  h\right)
}\right\Vert _{\infty}<M\quad\mathrm{or}\quad\limsup_{h\rightarrow0}\left\Vert
\frac{\Gamma_{K}}{v\left(  h\right)  }\right\Vert _{\infty}=+\infty.
\]

\end{itemize}

The following example will illustrate the points above. Take $T$ a diagonal
operator expressed in a basis of $H$ and defined this way : $T\left(
h\right)  =diag\left(  a_{i}\left(  h\right)  \right)  $ where%
\[
a_{i}\left(  h\right)  =h\frac{\lambda_{i}}{\lambda_{i}+h}+\frac{h^{3/2}%
}{\lambda_{i}+h},
\]
$h\downarrow0$ and $\lambda\in l_{1}$. The reader will be easily convinced
that we are in a situation similar to the one of theorem \ref{LCO} : the
$i^{th}$ cell of the bounded operator $T$ is asymptotically equivalent with
$h+h^{3/2}/\lambda_{i}.$ Here $v\left(  h\right)  =h,$ $w\left(  h\right)
=h^{3/2}$ and $\mathcal{R=}diag\left(  \lambda_{i}^{-1}\right)  $ is
unbounded. Then it is elementary algebra to prove that :%
\begin{align}
\left\Vert T\left(  h\right)  -hI\right\Vert _{\infty}  &  \leq h^{1/2}%
\rightarrow0,\nonumber\\
\left\Vert \frac{T\left(  h\right)  }{h}\right\Vert _{\infty}  &
=1+h^{-1/2}\rightarrow+\infty. \label{gorillaz}%
\end{align}
However if
\begin{align*}
a_{i}\left(  h\right)   &  =h\frac{\lambda_{i}}{\lambda_{i}+h}+\frac{h^{2}%
}{\lambda_{i}+h},\\
\left\Vert \frac{T\left(  h\right)  }{h}\right\Vert _{\infty}  &  =2.
\end{align*}

Let us focus again on the local covariance operator $\Gamma_{K}$. If one tries
to bound its norm a first attempt gives :%
\begin{align}
\left\Vert \Gamma_{K}\right\Vert _{\infty}  &  \leq\mathbb{E}\left\Vert
K\left(  \frac{\left\Vert X-x_{0}\right\Vert }{h}\right)  \left(  \left(
X-x_{0}\right)  \otimes\left(  X-x_{0}\right)  \right)  \right\Vert _{\infty
}\nonumber\\
&  =\mathbb{E}\left[  K\left(  \frac{\left\Vert X-x_{0}\right\Vert }%
{h}\right)  \left\Vert X-x_{0}\right\Vert ^{2}\right]  . \label{az}%
\end{align}
The next theorem assesses that under mild conditions that hold in the gaussian
framework this bound is not sharp. At this point we should turn back to
Theorem \ref{LCO}, especially to the sequences $v\left(  h\right)  $ and
$w\left(  h\right)  $ mentioned within this Theorem. The next Proposition
provides first a bound then under an additional assumption an equivalent
sequence for $v\left(  h\right)  .$ The case of $w\left(  h\right)  $ will not
be treated since the inspection of the method of proof would easily lead to
similar results.

\begin{proposition}
\label{Compay}Let as above%
\[
v\left(  h\right)  =\mathbb{E}\left(  K\left(  \frac{\left\Vert X-x_{0}%
\right\Vert }{h}\right)  \left\Vert X-x_{0}\right\Vert \rho\left(  \left\Vert
X-x_{0}\right\Vert \right)  \right)
\]
then%
\[
\frac{v\left(  h\right)  }{\mathbb{E}\left[  K\left(  \frac{\left\Vert
X-x_{0}\right\Vert }{h}\right)  \left\Vert X-x_{0}\right\Vert ^{2}\right]
}\rightarrow0.
\]
Besides if $\rho$ is regularly varying at $0$ with positive index%
\[
v\left(  h\right)  \sim K\left(  1\right)  h\rho\left(  h\right)  F\left(
h\right)  .
\]

\end{proposition}

\begin{remark}
The auxiliary functions appearing within Proposition \ref{MED} are both
regularly varying with indices~:%
\[
d=\left(  3+4\alpha\right)  /\left(  1+2\alpha\right)
\]
for the first and%
\[
d=1
\]
for the second.
\end{remark}

It is time for us to state the second main result of this paper. This Theorem
is in a way complementary to the previous one. It provides an asymptotic first
order development of $\Gamma_{K}$.

\begin{theorem}
\label{KT}Suppose $\mathbf{A}_{1}$ and $\mathbf{A}_{\mathbf{2}}$ hold. Let
$\mathcal{V}_{0}$ be a fixed neighborhood of $0$ and%
\[
a_{i}=\sup_{t\in\mathcal{V}_{0}}\left\vert \frac{f_{i}\left(  t+\left\langle
x_{0},e_{i}\right\rangle \right)  -f_{i}\left(  \left\langle x_{0}%
,e_{i}\right\rangle \right)  }{f_{i}\left(  \left\langle x_{0},e_{i}%
\right\rangle \right)  }\right\vert .
\]
If
\begin{equation}
\sum_{i=1}^{+\infty}a_{i}^{2}<+\infty\label{mig}%
\end{equation}
we have :%
\[
\left\Vert \Gamma_{K}-v\left(  h\right)  I\right\Vert _{\infty}=O\left(
v\left(  h\right)  \right)  .
\]

\end{theorem}

\begin{remark}
By Proposition \ref{Compay}, the Theorem just provides a bound sharper than
the rather "na\"{\i}ve" one at display (\ref{az}) since obviously due to
(\ref{F2}). A better result would be to obtain a second order term which would
mean here to provide the explicit operator hidden behind the "$O\left(
v\left(  h\right)  \right)  $". This has still to be done and holds perhaps
under reinforced assumptions on both $x_{0}$ and the $f_{i}$'s. But Theorem
\ref{LCO}, as well as the examples treated above (see display (\ref{gorillaz}%
)) let us claim that : "if ever $\left[  \Gamma_{K}-v\left(  h\right)
I\right]  /s\left(  h\right)  $ converges in norm, then necessarily $s\left(
h\right)  =w\left(  h\right)  $ and the limiting operator is then
$\mathcal{R}$".
\end{remark}

\begin{remark}
It has been seen above that we could not expect to obtain a $O\left(  w\left(
h\right)  \right)  $ on the right, instead of $O\left(  v\left(  h\right)
\right)  $ because in many situations the operator $\mathcal{R}$ will be
unbounded. We see that the price to pay to enhance a "weak" result such as
Theorem \ref{LCO} to a "uniform" one such as Theorem \ref{KT} is a slower rate
of decrease since obviously%
\[
\frac{w\left(  h\right)  }{v\left(  h\right)  }\rightarrow0.
\]

\end{remark}

Assumption (\ref{mig}) must be commented and illustrated by investigating some examples.

\begin{example}
(Gauss) If $X$ is gaussian, straightforward computations lead to :%
\[
\frac{f_{i}\left(  t+\left\langle x_{0},e_{i}\right\rangle \right)
-f_{i}\left(  \left\langle x_{0},e_{i}\right\rangle \right)  }{f_{i}\left(
\left\langle x_{0},e_{i}\right\rangle \right)  }=\exp\left(  \frac
{-t^{2}-2\left\langle x_{0},e_{i}\right\rangle t}{2\lambda_{i}}\right)  -1
\]
and if $i$ is large enough (hence $\left\langle x_{0},e_{i}\right\rangle $
small enough),%
\[
a_{i}=\left\vert \exp\left(  \frac{\left\langle x_{0},e_{i}\right\rangle ^{2}%
}{2\lambda_{i}}\right)  -1\right\vert .
\]
Then if $\frac{\left\langle x_{0},e_{i}\right\rangle ^{2}}{\lambda_{i}}$ tends
to zero,%
\[
a_{i}\leq\frac{\left\langle x_{0},e_{i}\right\rangle ^{2}}{4\lambda_{i}}%
\]
and (\ref{mig}) holds when%
\begin{equation}
\sum_{i=1}^{+\infty}\frac{\left\langle x_{0},e_{i}\right\rangle ^{4}}%
{\lambda_{i}^{2}}<+\infty. \label{Hotel}%
\end{equation}

\end{example}

\begin{example}
(Laplace) If%
\[
f_{i}\left(  t\right)  =\frac{1}{2\lambda_{i}}\exp\left(  -\frac{\left\vert
t\right\vert }{\lambda_{i}}\right)
\]
we have
\[
\frac{f_{i}\left(  t+\left\langle x_{0},e_{i}\right\rangle \right)
-f_{i}\left(  \left\langle x_{0},e_{i}\right\rangle \right)  }{f_{i}\left(
\left\langle x_{0},e_{i}\right\rangle \right)  }=\exp\left(  \frac{\left\vert
\left\langle x_{0},e_{i}\right\rangle \right\vert -\left\vert \left\langle
x_{0},e_{i}\right\rangle -t\right\vert }{2\lambda_{i}}\right)  -1
\]
and assumption (\ref{mig}) holds when%
\[
\sum_{i=1}^{+\infty}\left[  \exp\left(  \frac{\left\vert \left\langle
x_{0},e_{i}\right\rangle \right\vert }{2\lambda_{i}}\right)  -1\right]
^{2}<+\infty,
\]
hence when%
\[
\sum_{i=1}^{+\infty}\frac{\left\langle x_{0},e_{i}\right\rangle ^{2}}%
{\lambda_{i}^{2}}<+\infty.
\]

\end{example}

\begin{example}
(Unimodal densities) Since $X$ is assumed to be centered and the $f_{i}$'s are
the densities of the random variables $\left\langle X-x_{0},e_{i}\right\rangle
,$ we may extend both previous examples to a slightly more general situation.
Indeed both variance and expectation of $\left\langle X-x_{0},e_{i}%
\right\rangle $ tend to zero (they are respectively $\lambda_{i}$ and
$-\left\langle x_{0},e_{i}\right\rangle $) and we can consider the case when
$f_{i}$ features a single peak (a mode) at $-\left\langle x_{0},e_{i}%
\right\rangle $) and concentrates around $0.$ Then%
\[
a_{i}=\left\vert \frac{f_{i}\left(  0\right)  -f_{i}\left(  \left\langle
x_{0},e_{i}\right\rangle \right)  }{f_{i}\left(  \left\langle x_{0}%
,e_{i}\right\rangle \right)  }\right\vert .
\]
If we try to go beyond this relationship, a simple development of $f_{i}$
around zero provides :%
\[
a_{i}=-\left\langle x_{0},e_{i}\right\rangle \frac{f_{i}^{\prime}\left(
\left\langle x_{0},e_{i}\right\rangle \right)  }{f_{i}\left(  \left\langle
x_{0},e_{i}\right\rangle \right)  }+\frac{\left\langle x_{0},e_{i}%
\right\rangle ^{2}}{2}\frac{f_{i}^{\prime\prime}\left(  c_{i}\right)  }%
{f_{i}\left(  \left\langle x_{0},e_{i}\right\rangle \right)  }%
\]
where $c_{i}$ lies somewhere between $0$ and $\left\langle x_{0}%
,e_{i}\right\rangle $ Then it is plain to see that assumption (\ref{mig})
turns out to hold when : on a first hand (applying Cauchy-Schwartz inequality)
condition $\mathbf{C}_{2}$ holds and on the other hand when%
\[
\sup_{i}\frac{f_{i}^{\prime\prime}\left(  c_{i}\right)  }{f_{i}\left(
\left\langle x_{0},e_{i}\right\rangle \right)  }<+\infty
\]
which is not exactly assumption $\mathbf{C}_{1}^{\prime}$ but which is not
that far. Developping $f_{i}$ up to$f_{i}^{\prime\prime\prime}$ would let
$\mathbf{C}_{1}^{\prime}$ appear but would also creates an additional term.
This does not prove that (\ref{mig}) is a necessary and sufficient condition
for convergence of $\Gamma_{K}$ in the sense of Theorem \ref{KT}. But the
closeness of (\ref{mig}) with the "weak" conditions $\mathbf{C}_{1}^{\prime}$
and $\mathbf{C}_{2}$ shows that Theorem \ref{KT} is obtained under rather mild assumptions.
\end{example}

\begin{remark}
It has been assumed throughout the paper that $X$ has null expectation. When
$\mathbb{E}X\neq0,$ the theorems continue to hold but assumptions such as
(\ref{mig}) may implicitely involve this expectation itself. For instance
(\ref{Hotel}) is replaced by :%
\[
\sum_{i=1}^{+\infty}\frac{\left\langle x_{0}+\mu,e_{i}\right\rangle ^{4}%
}{\lambda_{i}^{2}}<+\infty
\]
where $\mu=\mathbb{E}X$. More generally speaking in the situation when the
$X^{\prime}s$ are not centered, a condition involving $x_{0}$ should be
replaced by the same condition involving $x_{0}+\mu.$
\end{remark}

\subsection{Convergence of the empirical covariance operator}

We go on with asymptotics for the empirical covariance operators, namely mean
square error. This part is short since the most delicate issue was sorted out
in the previous section.

\begin{theorem}
\label{covop}When assumptions $\mathbf{A}_{1}$and $\mathbf{A}_{2}$ hold the
following asymptotic results are true :%
\[
\mathbb{E}\left\Vert \Gamma_{K,n}-\Gamma_{K}\right\Vert _{\infty}^{2}%
=\frac{h^{4}F\left(  h\right)  }{n}K^{2}\left(  1\right)  \left(  1+o\left(
1\right)  \right)
\]

\end{theorem}

The same sort of results hold for the eigenelements. Recall that $\left(
\lambda_{n,p}\right)  _{p\in\mathbb{N}}$ (resp. $\left(  \lambda_{p}\right)
_{p\in\mathbb{N}}$) stands for the eigenvalues of $\Gamma_{K,n}$ (resp.
$\Gamma_{K}$) and $\left(  \pi_{n,p}\right)  _{p\in\mathbb{N}}$ (resp.
$\left(  \pi_{p}\right)  _{p\in\mathbb{N}}$) stands for the associated
eigenprojectors. The following Theorem estimates the rate of decrease to zero
for the eigenelements.

\begin{corollary}
\label{spec}Under the same conditions as in Theorem \ref{covop}, for fixed
$p\in\mathbb{N}$,%
\begin{align*}
\mathbb{E}\left(  \lambda_{n,p}-\lambda_{p}\right)  ^{2}  &  =O\left(
\frac{h^{4}F\left(  h\right)  }{n}\right) \\
\mathbb{E}\left\Vert \pi_{n,p}-\pi_{p}\right\Vert _{\infty}^{2}  &  =O\left(
\frac{h^{4}F\left(  h\right)  }{n}\right)
\end{align*}

\end{corollary}

Only a sketch of the proof of this Corollary is given since it may be seen as
a by-product of an article by Mas and Menneteau (2003). Under simple
additional assumptions, exact constants could be computed in both above
displays by applying the formulas that appear in Theorem 1.2 p.129 in their
article. But these computations are beyond the scope of this article : they
make it necessary to introduce and explain Kato's perturbation theory as well
as the associated functional calculus for linear operators. The interested
reader is referred to Kato (1976), Dunford-Schwartz (1988), Gohberg, Goldberg
and Kaashoek (1991).

\subsection{Conclusion and perspectives}

We proposed a first approach to what we named "local covariance operators".
This article does not aim at giving an exhaustive list of their features but
provides some clues for further and deeper study. For instance we can take for
granted that small probabilities naturally appear when estimating the rates of
convergence and that the class of $\Gamma$-varying functions ( at $0$)
provides an accurate setting. It also turns out from Theorem \ref{LCO} that
the asymptotic behaviour of $\Gamma_{K}$ is quite unusual and let appear,
through $\mathcal{R}$, an unbounded operator operator whereas $\Gamma_{K}$
tend s to zero.

Several issues will have to be addressed in the future. We can list some of
them. The case of first order truncated moments, that is :%
\[
\mathbb{E}\left[  \left(  X-x_{0}\right)  K\left(  \left\Vert \frac{X-x_{0}%
}{h}\right\Vert \right)  \right]  ,
\]
could certainly be studied in the framework of this article and with the same
computational techniques. Almost sure convergence as well as weak convergence
(convergence in distribution) could be adressed by following the same lines.

Another crucial issue is the existence and the properties of the inverse of
$\Gamma_{K}$ when it exists. Indeed let us introduce the nonparametric
regression model for functional random variables :%
\[
y=r\left(  X\right)  +\varepsilon,
\]
where $\left(  y,X\right)  \in\mathbb{R}\times H.$ Investigating a pointwise
estimate $\widehat{r}\left(  x_{0}\right)  $ of $r\left(  x_{0}\right)  $ by
local linear methods leads to finding the inverse (or a pseudo-inverse) of
$\Gamma_{K}.$

\section{Mathematical derivations}

For any $x=\sum x_{k}e_{k}$ in $H$ and for $\left(  i,j\right)  \in
\mathbb{N}^{2}$, $i\neq j,$ set%
\begin{align*}
\left\Vert x\right\Vert _{\neq i}^{2}  &  =\sum_{k\neq i}x_{k}^{2}\\
\left\Vert x\right\Vert _{\neq ij}^{2}  &  =\sum_{k\neq i,j}x_{k}^{2}%
\end{align*}
and denote $f_{\neq i}$ the density of $\left\Vert X\right\Vert _{\neq i}$ as
well as $f_{\neq ij}$ the density of $\left\Vert X\right\Vert _{\neq ij}$. It
is clear that when assumption $\mathbf{A}_{1}$ holds $\left\langle
X,e_{i}\right\rangle $ and $\left\Vert X\right\Vert _{\neq i}$ are independent
random variables.

In order to alleviate the notations, within the proofs -unless explicitely
mentioned- $X$ will stand for $X-x_{0}$ ($x_{0}$ is dropped since it is fixed
but we keep aware that all our results and notations, especially small ball
probabilities depend on $x_{0}$) and $f_{i}$ for the density of $\left\langle
X-x_{0},e_{i}\right\rangle $.

\subsection{Preliminary material}

We begin with preliminary Lemmas which are assessed in a general setting and
will be applied later. We recall the definition of the Gamma function $\Gamma$
:%
\[
\Gamma\left(  u\right)  =\int_{0}^{+\infty}s^{u-1}\exp\left(  -s\right)  ds.
\]

\begin{lemma}
\label{prel}If $f$ belongs to the class $\Gamma$ with auxiliary function
$\rho$, then for all $p\in\mathbb{N}$,%
\[
\int_{0}^{1}\frac{t^{p}}{\sqrt{1-t^{2}}}f\left(  s\sqrt{1-t^{2}}\right)
dt\underset{s\rightarrow0}{\sim}2^{\frac{p-1}{2}}\Gamma\left(  \frac{p+1}%
{2}\right)  f\left(  s\right)  \left(  \frac{\rho\left(  s\right)  }%
{s}\right)  ^{\frac{p+1}{2}}.
\]

\end{lemma}

\begin{proof}
We start with the following change of variable : $s\sqrt{1-t^{2}}%
=s-\rho\left(  s\right)  x$%
\begin{align*}
&  \int_{0}^{1}\frac{t^{p}}{\sqrt{1-t^{2}}}f\left(  s\sqrt{1-t^{2}}\right)
dt\\
&  =\int_{0}^{s/\rho\left(  s\right)  }\frac{\rho\left(  s\right)  }{s}\left(
\frac{\rho\left(  s\right)  }{s}x\left(  2+\frac{\rho\left(  s\right)  }%
{s}\right)  \right)  ^{\frac{p-1}{2}}f\left(  s-\rho\left(  s\right)
x\right)  dx\\
&  =2^{\frac{p-1}{2}}\left(  \frac{\rho\left(  s\right)  }{s}\right)
^{\frac{p+1}{2}}f\left(  s\right)  \int_{0}^{s/\rho\left(  s\right)  }\left(
x\left(  1+\frac{\rho\left(  s\right)  }{2s}\right)  \right)  ^{\frac{p-1}{2}%
}\frac{f\left(  s-\rho\left(  s\right)  x\right)  }{f\left(  s\right)
\exp\left(  -x\right)  }\exp\left(  -x\right)  dx.
\end{align*}
To conclude it suffices to prove that
\[
\frac{\int_{0}^{s/\rho\left(  s\right)  }\left(  x\left(  1+\frac{\rho\left(
s\right)  }{2s}\right)  \right)  ^{\frac{p-1}{2}}\frac{f\left(  s-\rho\left(
s\right)  x\right)  }{f\left(  s\right)  \exp\left(  -x\right)  }\exp\left(
-x\right)  dx-\Gamma\left(  \frac{p+1}{2}\right)  }{\Gamma\left(  \frac
{p+1}{2}\right)  }\underset{s\rightarrow0}{\rightarrow}0.
\]
But the numerator may be rewritten%
\begin{equation}
\int_{0}^{s/\rho\left(  s\right)  }\left[  \left(  \left(  1+\frac{\rho\left(
s\right)  }{2s}\right)  \right)  ^{\frac{p-1}{2}}\frac{f\left(  s-\rho\left(
s\right)  x\right)  }{f\left(  s\right)  \exp\left(  -x\right)  }-1\right]
x^{\frac{p-1}{2}}\exp\left(  -x\right)  dx. \label{marine}%
\end{equation}
Let us study briefly the sequence of functions $\frac{f\left(  s-\rho\left(
s\right)  x\right)  }{f\left(  s\right)  \exp\left(  -x\right)  }.$ The sign
of the first order derivative is the sign of :%
\[
\frac{f\left(  s-\rho\left(  s\right)  x\right)  }{f^{\prime}\left(
s-\rho\left(  s\right)  x\right)  }-\rho\left(  s\right)
\]
By Theorem 3.10.11 in Bingham, Goldie, Teugels (1987), $f^{\prime}$ is
$\Gamma$-varying with same auxiliary function as $f$ and by Corollary 3.10.5
(b) p.177 ibidem we know that%
\[
\rho\left(  s-\rho\left(  s\right)  x\right)  =\frac{f\left(  s-\rho\left(
s\right)  x\right)  }{f^{\prime}\left(  s-\rho\left(  s\right)  x\right)  }%
\]
Then since $\rho\left(  0\right)  =0$ and $\rho\geq0,$ $\rho$ is strictly
increasing in a neighborhood of zero and $\rho\left(  s-\rho\left(  s\right)
x\right)  \leq\rho\left(  s\right)  ,$ $\frac{f\left(  s-\rho\left(  s\right)
x\right)  }{f\left(  s\right)  \exp\left(  -x\right)  }$ is nonincreasing on
$\left[  0,s/\rho\left(  s\right)  \right]  $ (as a function of $x$) hence%
\[
\sup_{x\in\left[  0,s/\rho\left(  s\right)  \right]  }\frac{f\left(
s-\rho\left(  s\right)  x\right)  }{f\left(  s\right)  \exp\left(  -x\right)
}=1
\]

Now together with display (\ref{defgamma}) and (\ref{F2}) we can apply
Lebesgue's dominated convergence Theorem to which completes the proof of
(\ref{marine}) the Lemma.
\end{proof}

The next lemma

If $U$ and $V$are two real valued random variables, $f_{U,V}$ denotes the
joint density of the couple $\left(  U,V\right)  $. We need to compute four
densities :

\begin{lemma}
\label{dens}We have :%
\begin{align}
f_{\left\langle X,e_{i}\right\rangle ,\left\Vert X\right\Vert }\left(
u,v\right)   &  =\frac{v}{\sqrt{v^{2}-u^{2}}}f_{i}\left(  u\right)  f_{\neq
i}\left(  \sqrt{v^{2}-u^{2}}\right)  1\hspace{-3pt}1_{\left\{  v\geq\left\vert
u\right\vert \right\}  },\label{d1}\\
f_{\left\Vert X\right\Vert }\left(  v\right)   &  =v\int_{-1}^{1}\frac
{f_{i}\left(  vt\right)  }{\sqrt{1-t^{2}}}f_{\neq i}\left(  v\sqrt{1-t^{2}%
}\right)  dt \label{d1bis}%
\end{align}
and%
\begin{align}
&  f_{\left\langle X,e_{i}\right\rangle ,\left\langle X,e_{j}\right\rangle
,\left\Vert X\right\Vert }\left(  t,u,v\right) \nonumber\\
&  =\frac{v}{\sqrt{v^{2}-u^{2}-t^{2}}}f_{i}\left(  t\right)  f_{j}\left(
u\right)  f_{\neq ij}\left(  \sqrt{v^{2}-u^{2}-t^{2}}\right)  1\hspace
{-3pt}1_{\left\{  v\geq\sqrt{u^{2}+t^{2}}\right\}  },\label{d3}\\
&  f_{\left\Vert X\right\Vert }\left(  v\right) \nonumber\\
&  =v^{2}\int_{0}^{2\pi}\int_{0}^{1}\frac{x}{\sqrt{1-x^{2}}}f_{i}\left(
vx\cos\theta\right)  f_{j}\left(  vx\sin\theta\right)  f_{\neq ij}\left(
v\sqrt{1-x^{2}}\right)  dxd\theta.\nonumber
\end{align}

\end{lemma}

\begin{proof}
We only compute the first density since we could get the second by
integration. The third and fourth could be obtained the same way. Let $h$ be
any bounded measurable function.%
\begin{align*}
\mathbb{E}h\left(  \left\langle X,e_{i}\right\rangle ,\left\Vert X\right\Vert
\right)   &  =\int h\left(  x_{i},\sqrt{x_{i}^{2}+y_{i}^{2}}\right)
f_{i}\left(  x_{i}\right)  f_{\neq i}\left(  y_{i}\right)  dx_{i}dy_{i}\\
&  =\int v\frac{h\left(  u,v\right)  }{\sqrt{v^{2}-u^{2}}}f_{i}\left(
u\right)  f_{\neq i}\left(  \sqrt{v^{2}-u^{2}}\right)  1\hspace{-3pt}%
1_{\left\{  v\geq\left\vert u\right\vert \right\}  }dudv\\
&  =\int h\left(  u,v\right)  f_{\left\langle X,e_{i}\right\rangle ,\left\Vert
X\right\Vert }\left(  u,v\right)  dudv.
\end{align*}
Identifying both last terms we get%
\[
f_{\left\langle X,e_{i}\right\rangle ,\left\Vert X\right\Vert }\left(
u,v\right)  =\frac{v}{\sqrt{v^{2}-u^{2}}}f_{i}\left(  u\right)  f_{\neq
i}\left(  \sqrt{v^{2}-u^{2}}\right)  1\hspace{-3pt}1_{\left\{  v\geq\left\vert
u\right\vert \right\}  }.
\]
Integrating this density with respect to the variable $u$ yields
$f_{\left\Vert X\right\Vert }\left(  v\right)  $ as in (\ref{d1bis}).
\end{proof}

\begin{lemma}
\label{Loundge}The following hold :%
\begin{align*}
&  f_{\left\Vert X\right\Vert }\left(  v\right)  \underset{0}{\sim}%
\Gamma\left(  \frac{1}{2}\right)  \sqrt{2v\rho\left(  v\right)  }f_{i}\left(
0\right)  f_{\neq i}\left(  v\right)  ,\\
&  f_{\left\Vert X\right\Vert }\left(  v\right)  \underset{0}{\sim}2\pi
f_{i}\left(  0\right)  f_{j}\left(  0\right)  v\rho\left(  v\right)  f_{\neq
ij}\left(  v\right)  .
\end{align*}
Besides if $f_{\left\Vert X\right\Vert }$, $f_{\neq i}$ and $f_{\neq ij}$ are
$\Gamma$-varying for all $i$ and $j$ then they have all $\rho$ as auxiliary function.
\end{lemma}

\begin{proof}
We restrict to proving the Lemma for $f_{\neq i}.$ From Lemma \ref{prel} and
(\ref{d1bis}) we get%
\begin{align*}
f_{\left\Vert X\right\Vert }\left(  v\right)   &  =v\int_{-1}^{1}\frac
{f_{i}\left(  vt\right)  }{\sqrt{1-t^{2}}}f_{\neq i}\left(  v\sqrt{1-t^{2}%
}\right)  dt\\
&  \underset{0}{\sim}2vf_{i}\left(  0\right)  \int_{0}^{1}\frac{1}%
{\sqrt{1-t^{2}}}f_{\neq i}\left(  v\sqrt{1-t^{2}}\right)  dt\\
&  \underset{0}{\sim}\Gamma\left(  \frac{1}{2}\right)  \sqrt{2v\rho_{i}\left(
v\right)  }f_{i}\left(  0\right)  f_{\neq i}\left(  v\right)
\end{align*}
where $\rho_{i}$ denotes the auxiliary function of $f_{\neq i}$. Now we also
have for all $x\geq0$ :%
\begin{align*}
&  \frac{f_{\left\Vert X\right\Vert }\left(  v+\rho_{i}\left(  v\right)
x\right)  }{f_{\left\Vert X\right\Vert }\left(  v\right)  }\underset{0}{\sim
}\frac{\sqrt{2\left(  v+\rho_{i}\left(  v\right)  x\right)  \rho_{i}\left(
v+\rho_{i}\left(  v\right)  x\right)  }f_{i}\left(  0\right)  f_{\neq
i}\left(  v+\rho_{i}\left(  v\right)  x\right)  }{\sqrt{2v\rho_{i}\left(
v\right)  }f_{i}\left(  0\right)  f_{\neq i}\left(  v\right)  }\\
&  \underset{0}{\sim}\sqrt{\frac{\left(  v+\rho_{i}\left(  v\right)  x\right)
\rho_{i}\left(  v+\rho_{i}\left(  v\right)  x\right)  }{v\rho_{i}\left(
v\right)  }}\exp\left(  x\right)  .
\end{align*}
The term%
\[
\frac{\left(  v+\rho_{i}\left(  v\right)  x\right)  \rho_{i}\left(  v+\rho
_{i}\left(  v\right)  x\right)  }{v\rho_{i}\left(  v\right)  }=\left(
1+\frac{\rho_{i}\left(  v\right)  x}{v}\right)  \left(  \frac{\rho_{i}\left(
v+\rho_{i}\left(  v\right)  x\right)  }{\rho_{i}\left(  v\right)  }\right)
\]
tends to $1$ by \textbf{Fact 2}. Finally%
\[
\frac{f_{\left\Vert X\right\Vert }\left(  v+\rho_{i}\left(  v\right)
x\right)  }{f_{\left\Vert X\right\Vert }\left(  v\right)  }\rightarrow
\exp\left(  x\right)
\]
and $\rho_{i}$ is also the auxiliary function for $f_{\left\Vert X\right\Vert
}$. Since $f_{\left\Vert X\right\Vert }$ is also $\Gamma$-varying with
auxiliary function $\rho$, we can set $\rho_{i}=\rho$ (the auxiliary function
is unique up to an asymptotic equivalence, see Corollary 3.10.5 (b) p.177 in
Bingham, Goldie, Teugels (1987)). The same steps would lead us to the second
part of the Lemma.
\end{proof}

\subsection{Proof of cell-by-cell results}

By "cell" we just mean that, identifiying $\Gamma_{K}$ with an infinite
matrix, we consider in this subsection asymptotics for $\left\langle
\Gamma_{K}e_{i},e_{j}\right\rangle .$ We study first the diagonal of
$\Gamma_{K}:$

\begin{proposition}
\label{T1}Fix the index $i\in\mathbb{N}$ :%
\begin{align}
&  \mathbb{E}\left[  K\left(  \frac{\left\Vert X\right\Vert }{h}\right)
\left\langle X,e_{i}\right\rangle ^{2}\right] \nonumber\\
&  \underset{h\rightarrow0}{\sim}\int_{0}^{h}v\rho\left(  v\right)  K\left(
\frac{v}{h}\right)  f_{\left\Vert X\right\Vert }\left(  v\right)
dv+\frac{f_{i}^{\prime\prime}\left(  0\right)  }{f_{i}\left(  0\right)  }%
\int_{0}^{h}v^{2}\rho^{2}\left(  v\right)  K\left(  \frac{v}{h}\right)
f_{\left\Vert X\right\Vert }\left(  v\right)  dv. \label{indetendances}%
\end{align}
\bigskip
\end{proposition}

\begin{corollary}
From the above we deduce that :%
\[
\mathbb{E}\left[  \left\langle X,e_{i}\right\rangle ^{2}|\left\Vert
X\right\Vert =v\right]  \sim v\rho\left(  v\right)  .
\]
\bigskip
\end{corollary}

\begin{proof}
We start from the joint density at display (\ref{d1}) :%

\begin{align}
&  \mathbb{E}\left[  K\left(  \frac{\left\Vert X\right\Vert }{h}\right)
\left\langle X,e_{i}\right\rangle ^{2}\right] \nonumber\\
&  =\int\int K\left(  \frac{v}{h}\right)  u^{2}\frac{v}{\sqrt{v^{2}-u^{2}}%
}f_{i}\left(  u\right)  f_{\neq i}\left(  \sqrt{v^{2}-u^{2}}\right)
1\hspace{-3pt}1_{\left\{  h\geq v\geq\left\vert u\right\vert \right\}
}dudv\nonumber\\
&  =\int_{0}^{h}vK\left(  \frac{v}{h}\right)  \left(  \int_{-v}^{v}\frac
{u^{2}}{\sqrt{v^{2}-u^{2}}}f_{i}\left(  u\right)  f_{\neq i}\left(
\sqrt{v^{2}-u^{2}}\right)  du\right)  dv\label{Malcolm}\\
&  =\int_{0}^{h}v^{3}K\left(  \frac{v}{h}\right)  \left(  \int_{-1}^{1}%
\frac{x^{2}}{\sqrt{1-x^{2}}}f_{i}\left(  xv\right)  f_{\neq i}\left(
v\sqrt{1-x^{2}}\right)  dx\right)  dv.\nonumber
\end{align}
Setting :%
\[
I_{i}\left(  v\right)  =\frac{\int_{-1}^{1}\frac{vx^{2}}{\sqrt{1-x^{2}}}%
f_{i}\left(  xv\right)  f_{\neq i}\left(  v\sqrt{1-x^{2}}\right)  dx}%
{\int_{-1}^{1}\frac{v}{\sqrt{1-x^{2}}}f_{i}\left(  vx\right)  f_{\neq
i}\left(  v\sqrt{1-x^{2}}\right)  dx},
\]
we have%
\begin{equation}
\mathbb{E}\left[  K\left(  \frac{\left\Vert X\right\Vert }{h}\right)
\left\langle X,e_{i}\right\rangle ^{2}\right]  =\int_{0}^{h}v^{2}K\left(
\frac{v}{h}\right)  I_{i}\left(  v\right)  f_{\left\Vert X\right\Vert }\left(
v\right)  dv. \label{inde1}%
\end{equation}
We focus on $I\left(  v\right)  $ and prove that :%
\begin{equation}
I_{i}\left(  v\right)  \sim\frac{\rho\left(  v\right)  }{v}\left(
1+v\rho\left(  v\right)  \frac{f_{i}^{\prime\prime}\left(  0\right)  }%
{f_{i}\left(  0\right)  }\right)  . \label{inde2}%
\end{equation}
It is clear that from (\ref{inde1}) and (\ref{inde2}) we can derive
(\ref{indetendances}). \newline Now from $f_{i}\left(  xv\right)
=f_{i}\left(  0\right)  +vxf_{i}^{\prime}\left(  0\right)  +\frac{v^{2}x^{2}%
}{2}f_{i}^{\prime\prime}\left(  0\right)  +o\left(  v^{2}\right)  ,$ we can
develop $I_{i}\left(  v\right)  $. Setting :%
\[
\mathcal{J}_{p}\left(  v\right)  =\left(  \int_{-1}^{1}\frac{x^{p}}%
{\sqrt{1-x^{2}}}f_{\neq i}\left(  v\sqrt{1-x^{2}}\right)  dx\right)  \quad
p\in\mathbb{N},
\]
we get :%
\begin{align*}
I_{i}\left(  v\right)   &  =\frac{f_{i}\left(  0\right)  \mathcal{J}%
_{2}\left(  v\right)  +\frac{v^{2}}{2}f_{i}^{\prime\prime}\left(  0\right)
\mathcal{J}_{4}\left(  v\right)  }{f_{i}\left(  0\right)  \mathcal{J}%
_{0}\left(  v\right)  +\frac{v^{2}}{2}f_{i}^{\prime\prime}\left(  0\right)
\mathcal{J}_{2}\left(  v\right)  }\left(  1+o\left(  1\right)  \right) \\
&  =\frac{\mathcal{J}_{2}\left(  v\right)  }{\mathcal{J}_{0}\left(  v\right)
}\frac{1+\frac{v^{2}}{2}\frac{f_{i}^{\prime\prime}\left(  0\right)  }%
{f_{i}\left(  0\right)  }\frac{\mathcal{J}_{4}\left(  v\right)  }%
{\mathcal{J}_{2}\left(  v\right)  }}{1+\frac{v^{2}}{2}\frac{f_{i}%
^{\prime\prime}\left(  0\right)  }{f_{i}\left(  0\right)  }\frac
{\mathcal{J}_{2}\left(  v\right)  }{\mathcal{J}_{0}\left(  v\right)  }}\left(
1+o\left(  1\right)  \right) \\
&  =\frac{\mathcal{J}_{2}\left(  v\right)  }{\mathcal{J}_{0}\left(  v\right)
}\left(  1+\frac{v^{2}}{2}\frac{f_{i}^{\prime\prime}\left(  0\right)  }%
{f_{i}\left(  0\right)  }\left(  \frac{\mathcal{J}_{4}\left(  v\right)
}{\mathcal{J}_{2}\left(  v\right)  }-\frac{\mathcal{J}_{2}\left(  v\right)
}{\mathcal{J}_{0}\left(  v\right)  }\right)  \right)  \left(  1+o\left(
1\right)  \right)  .
\end{align*}
Then we invoke Lemma \ref{prel} and Lemma \ref{Loundge} to get :%
\begin{align*}
\frac{\mathcal{J}_{2}\left(  v\right)  }{\mathcal{J}_{0}\left(  v\right)  }
&  =\frac{\rho\left(  v\right)  }{v}\left(  1+o\left(  1\right)  \right)  ,\\
\frac{\mathcal{J}_{4}\left(  v\right)  }{\mathcal{J}_{2}\left(  v\right)  }
&  =3\frac{\rho\left(  v\right)  }{v}\left(  1+o\left(  1\right)  \right)  ,
\end{align*}
hence%
\[
I_{i}\left(  v\right)  =\frac{\rho\left(  v\right)  }{v}\left(  1+\frac
{f_{i}^{\prime\prime}\left(  0\right)  }{f_{i}\left(  0\right)  }v\rho\left(
v\right)  \right)  \left(  1+o\left(  1\right)  \right)
\]
which finishes the proof of the Proposition \ref{T1}.
\end{proof}

\begin{proposition}
\label{T2}Let us take $i\neq j$ in $\mathbb{N},$ we have :%
\[
\mathbb{E}\left[  K\left(  \frac{\left\Vert X\right\Vert }{h}\right)
\left\langle X,e_{i}\right\rangle \left\langle X,e_{j}\right\rangle \right]
\sim\frac{f_{i}^{\prime}\left(  0\right)  f_{j}^{\prime}\left(  0\right)
}{f_{i}\left(  0\right)  f_{j}\left(  0\right)  }\int_{0}^{h}K\left(  \frac
{v}{h}\right)  v^{2}\rho^{2}\left(  v\right)  f_{\left\Vert X\right\Vert
}\left(  v\right)  dv.
\]

\end{proposition}

\begin{proof}
The proof basically follows the same lines as the previous Proposition with a
few changes :%
\begin{subequations}
\begin{align}
&  \mathbb{E}\left[  K\left(  \frac{\left\Vert X\right\Vert }{h}\right)
\left\langle X,e_{i}\right\rangle \left\langle X,e_{j}\right\rangle \right]
\nonumber\\
&  =\int\int\int K\left(  \frac{v}{h}\right)  \frac{utv}{\sqrt{v^{2}%
-u^{2}-t^{2}}}f_{i}\left(  t\right)  f_{j}\left(  u\right)  f_{\neq ij}\left(
\sqrt{v^{2}-u^{2}-t^{2}}\right)  1\hspace{-3pt}1_{\left\{  h\geq v\geq
\sqrt{u^{2}+t^{2}}\right\}  }dudvdt\nonumber\\
&  =\int_{0}^{h}vK\left(  \frac{v}{h}\right)  \left(  \int\int\frac{ut}%
{\sqrt{v^{2}-u^{2}-t^{2}}}f_{i}\left(  t\right)  f_{j}\left(  u\right)
f_{\neq ij}\left(  \sqrt{v^{2}-u^{2}-t^{2}}\right)  1\hspace{-3pt}1_{\left\{
v\geq\sqrt{u^{2}+t^{2}}\right\}  }dudt\right)  dv\nonumber\\
&  =\int_{0}^{h}vK\left(  \frac{v}{h}\right) \nonumber\\
&  \times\left(  \int\int\frac{r^{3}\sin\theta\cos\theta}{\sqrt{v^{2}-r^{2}}%
}f_{i}\left(  r\cos\theta\right)  f_{j}\left(  r\sin\theta\right)  f_{\neq
ij}\left(  \sqrt{v^{2}-r^{2}}\right)  1\hspace{-3pt}1_{\left\{  v\geq
r\right\}  }drd\theta\right)  dv\label{th}\\
&  =\int_{0}^{h}v^{4}K\left(  \frac{v}{h}\right) \nonumber\\
&  \times\left(  \int\int\frac{x^{3}\sin\theta\cos\theta}{\sqrt{1-x^{2}}}%
f_{i}\left(  xv\cos\theta\right)  f_{j}\left(  xv\sin\theta\right)  f_{\neq
ij}\left(  v\sqrt{1-x^{2}}\right)  1\hspace{-3pt}1_{\left\{  1\geq
x\geq0\right\}  }dxd\theta\right)  dv\nonumber\\
&  =\int_{0}^{h}v^{2}K\left(  \frac{v}{h}\right)  J\left(  v\right)
f_{\left\Vert X\right\Vert }\left(  v\right)  dv, \label{bo}%
\end{align}
where%
\end{subequations}
\begin{equation}
J\left(  v\right)  =\frac{\int_{0}^{2\pi}\int_{0}^{1}\frac{x^{3}\sin\theta
\cos\theta}{\sqrt{1-x^{2}}}f_{i}\left(  xv\cos\theta\right)  f_{j}\left(
xv\sin\theta\right)  f_{\neq ij}\left(  v\sqrt{1-x^{2}}\right)  dxd\theta
}{\int_{0}^{2\pi}\int_{0}^{1}\frac{x}{\sqrt{1-x^{2}}}f_{i}\left(  vx\cos
\theta\right)  f_{j}\left(  vx\sin\theta\right)  f_{\neq ij}\left(
v\sqrt{1-x^{2}}\right)  dxd\theta}. \label{bo2}%
\end{equation}
At last we prove that :%
\[
J\left(  v\right)  \sim\frac{f_{i}^{\prime}\left(  0\right)  f_{j}^{\prime
}\left(  0\right)  }{f_{i}\left(  0\right)  f_{j}\left(  0\right)  }\rho
^{2}\left(  v\right)
\]
and go quickly through it :%
\begin{align*}
J\left(  v\right)   &  =\frac{\int_{0}^{2\pi}\int_{0}^{1}\frac{x^{3}\sin
\theta\cos\theta}{\sqrt{1-x^{2}}}f_{i}\left(  xv\cos\theta\right)
f_{j}\left(  xv\sin\theta\right)  f_{\neq ij}\left(  v\sqrt{1-x^{2}}\right)
1\hspace{-3pt}1_{\left\{  1\geq x\geq0\right\}  }dxd\theta}{\int_{0}^{2\pi
}\int_{0}^{1}\frac{x}{\sqrt{1-x^{2}}}f_{i}\left(  vx\cos\theta\right)
f_{j}\left(  vx\sin\theta\right)  f_{\neq ij}\left(  v\sqrt{1-x^{2}}\right)
dxd\theta}\\
&  =\frac{\int_{0}^{2\pi}\int_{0}^{1}\frac{x^{3}\sin\theta\cos\theta}%
{\sqrt{1-x^{2}}}f_{i}\left(  xv\cos\theta\right)  f_{j}\left(  xv\sin
\theta\right)  f_{\neq ij}\left(  v\sqrt{1-x^{2}}\right)  1\hspace
{-3pt}1_{\left\{  1\geq x\geq0\right\}  }dxd\theta}{f_{i}\left(  0\right)
f_{j}\left(  0\right)  \int_{0}^{2\pi}\int_{0}^{1}\frac{x}{\sqrt{1-x^{2}}%
}f_{\neq ij}\left(  v\sqrt{1-x^{2}}\right)  dxd\theta}\left(  1+o\left(
1\right)  \right) \\
&  =\frac{v^{2}}{2\pi}\frac{f_{i}^{\prime}\left(  0\right)  f_{j}^{\prime
}\left(  0\right)  }{f_{i}\left(  0\right)  f_{j}\left(  0\right)  }\frac
{\int_{0}^{2\pi}\sin^{2}\theta\cos^{2}\theta d\theta\int_{0}^{1}\frac{x^{5}%
}{\sqrt{1-x^{2}}}f_{\neq ij}\left(  v\sqrt{1-x^{2}}\right)  dx}{\int_{0}%
^{1}\frac{x}{\sqrt{1-x^{2}}}f_{\neq ij}\left(  v\sqrt{1-x^{2}}\right)
dx}\left(  1+o\left(  1\right)  \right) \\
&  =\frac{v^{2}}{2\pi}\frac{f_{i}^{\prime}\left(  0\right)  f_{j}^{\prime
}\left(  0\right)  }{f_{i}\left(  0\right)  f_{j}\left(  0\right)  }\frac{\pi
}{4}\frac{\mathcal{J}_{5}\left(  v\right)  }{\mathcal{J}_{1}\left(  v\right)
}\left(  1+o\left(  1\right)  \right) \\
&  =\frac{v^{2}}{8}\frac{f_{i}^{\prime}\left(  0\right)  f_{j}^{\prime}\left(
0\right)  }{f_{i}\left(  0\right)  f_{j}\left(  0\right)  }\frac{2^{2}%
\Gamma\left(  3\right)  \left(  \frac{\rho\left(  v\right)  }{v}\right)  ^{3}%
}{\Gamma\left(  1\right)  \left(  \frac{\rho\left(  v\right)  }{v}\right)
}\left(  1+o\left(  1\right)  \right) \\
&  =\frac{f_{i}^{\prime}\left(  0\right)  f_{j}^{\prime}\left(  0\right)
}{f_{i}\left(  0\right)  f_{j}\left(  0\right)  }\rho^{2}\left(  v\right)
\left(  1+o\left(  1\right)  \right)
\end{align*}
This last step ends the proof of Proposition \ref{T2}.\bigskip
\end{proof}

\textbf{Proof of Theorem \ref{LCO} :} The proof of the Theorem stems from
Propositions \ref{T1} and \ref{T2}. For instance,%
\[
\int_{0}^{h}v\rho\left(  v\right)  K\left(  \frac{v}{h}\right)  f_{\left\Vert
X\right\Vert }\left(  v\right)  dv=\mathbb{E}\left[  \left\Vert X\right\Vert
\rho\left(  \left\Vert X\right\Vert \right)  K\left(  \frac{\left\Vert
X\right\Vert }{h}\right)  \right]  .
\]

\subsection{Proof of norm results}

We decompose $\Gamma_{K}$ into two terms : a purely diagonal one and
non-diagonal one. In fact :%
\begin{align*}
\left\langle \Gamma_{K}^{d}e_{i},e_{i}\right\rangle  &  =\left\langle
\Gamma_{K}e_{i},e_{i}\right\rangle \\
\left\langle \Gamma_{K}^{d}e_{i},e_{j}\right\rangle  &  =0\quad i\neq j
\end{align*}
and
\begin{equation}
\Gamma_{K}^{\#d}=\Gamma_{K}-\Gamma_{K}^{d}. \label{di}%
\end{equation}
We first prove that :

\begin{lemma}
\label{chocolate}%
\[
\left\Vert \Gamma_{K}^{d}-v\left(  h\right)  I\right\Vert _{\infty}=O\left(
v\left(  h\right)  \right)
\]

\end{lemma}

It suffices to prove that :%
\[
\sup_{i\in\mathbb{N}}\left\vert \left\langle \Gamma_{K}e_{i},e_{i}%
\right\rangle -v\left(  h\right)  \right\vert =O\left(  v\left(  h\right)
\right)
\]
Let us denote :%
\[
\varphi_{i}\left(  t\right)  =\frac{f_{i}\left(  t\right)  -f_{i}\left(
0\right)  }{f_{i}\left(  0\right)  }%
\]
We have%
\[
\sup_{t\in\mathcal{V}_{0}}\left\vert \varphi_{i}\left(  t\right)  \right\vert
=a_{i}%
\]
where $a_{i}$ was introduced in Theorem \ref{KT}.

We start from (\ref{inde1}) and (\ref{inde2}) :%
\begin{align}
&  \left\langle \Gamma_{K}e_{i},e_{i}\right\rangle -v\left(  h\right)
\nonumber\\
&  =\mathbb{E}\left[  K\left(  \frac{\left\Vert X\right\Vert }{h}\right)
\left\langle X,e_{i}\right\rangle ^{2}\right]  -v\left(  h\right) \nonumber\\
&  =\int_{0}^{h}t^{2}K\left(  \frac{t}{h}\right)  I_{i}\left(  t\right)
f_{\left\Vert X\right\Vert }\left(  t\right)  dt-\int_{0}^{h}t\rho\left(
t\right)  K\left(  \frac{t}{h}\right)  f_{\left\Vert X\right\Vert }\left(
t\right)  dt\nonumber\\
&  =\int_{0}^{h}t^{2}K\left(  \frac{t}{h}\right)  \left(  I_{i}\left(
t\right)  -\frac{\rho\left(  t\right)  }{t}\right)  f_{\left\Vert X\right\Vert
}\left(  t\right)  dt \label{Sarko}%
\end{align}

We will first focus on :%
\[
I_{i}\left(  t\right)  -\frac{\rho\left(  t\right)  }{t}=I_{i}\left(
t\right)  -\frac{\mathcal{J}_{2}\left(  t\right)  }{\mathcal{J}_{0}\left(
t\right)  }+\frac{\mathcal{J}_{2}\left(  t\right)  }{\mathcal{J}_{0}\left(
t\right)  }-\frac{\rho\left(  t\right)  }{t}%
\]

Let us develop%
\begin{align*}
&  I_{i}\left(  t\right)  -\frac{\mathcal{J}_{2}\left(  t\right)
}{\mathcal{J}_{0}\left(  t\right)  }\\
&  =\frac{\int_{-1}^{1}\frac{x^{2}}{\sqrt{1-x^{2}}}f_{i}\left(  xt\right)
f_{\neq i}\left(  t\sqrt{1-x^{2}}\right)  dx}{\int_{-1}^{1}\frac{1}%
{\sqrt{1-x^{2}}}f_{i}\left(  tx\right)  f_{\neq i}\left(  t\sqrt{1-x^{2}%
}\right)  dx}-\frac{\int_{-1}^{1}\frac{x^{2}}{\sqrt{1-x^{2}}}f_{\neq i}\left(
t\sqrt{1-x^{2}}\right)  dx}{\int_{-1}^{1}\frac{1}{\sqrt{1-x^{2}}}f_{\neq
i}\left(  t\sqrt{1-x^{2}}\right)  dx}%
\end{align*}
It is plain that%
\[
I_{i}\left(  t\right)  =\frac{\int_{-1}^{1}\frac{x^{2}}{\sqrt{1-x^{2}}}\left(
1+\varphi_{i}\left(  tx\right)  \right)  f_{\neq i}\left(  t\sqrt{1-x^{2}%
}\right)  dx}{\int_{-1}^{1}\frac{1}{\sqrt{1-x^{2}}}\left(  1+\varphi
_{i}\left(  tx\right)  \right)  f_{\neq i}\left(  t\sqrt{1-x^{2}}\right)  dx}%
\]
Now denote%
\begin{align*}
\mathcal{J}_{0}^{\ast}\left(  t\right)   &  =\int_{-1}^{1}\frac{1}%
{\sqrt{1-x^{2}}}\varphi_{i}\left(  tx\right)  f_{\neq i}\left(  t\sqrt
{1-x^{2}}\right)  dx\\
\mathcal{J}_{2}^{\ast}\left(  t\right)   &  =\int_{-1}^{1}\frac{x^{2}}%
{\sqrt{1-x^{2}}}\varphi_{i}\left(  tx\right)  f_{\neq i}\left(  t\sqrt
{1-x^{2}}\right)  dx
\end{align*}
then%
\[
I_{i}\left(  t\right)  =\frac{\mathcal{J}_{0}\left(  t\right)  +\mathcal{J}%
_{0}^{\ast}\left(  t\right)  }{\mathcal{J}_{2}\left(  t\right)  +\mathcal{J}%
_{2}^{\ast}\left(  t\right)  }%
\]
hence%
\[
I_{i}\left(  t\right)  -\frac{\mathcal{J}_{2}\left(  t\right)  }%
{\mathcal{J}_{0}\left(  t\right)  }=\frac{\mathcal{J}_{2}^{\ast}\left(
t\right)  \mathcal{J}_{0}\left(  t\right)  -\mathcal{J}_{2}\left(  t\right)
\mathcal{J}_{0}^{\ast}\left(  t\right)  }{\mathcal{J}_{0}\left(  t\right)
\left(  \mathcal{J}_{0}\left(  t\right)  +\mathcal{J}_{0}^{\ast}\left(
t\right)  \right)  }%
\]
We are going to use the following inequalities :%
\begin{align*}
\left\vert \mathcal{J}_{0}^{\ast}\left(  t\right)  \right\vert  &  \leq
a_{i}\mathcal{J}_{0}\left(  t\right) \\
\left\vert \mathcal{J}_{2}^{\ast}\left(  t\right)  \right\vert  &  \leq
a_{i}\mathcal{J}_{2}\left(  t\right) \\
\mathcal{J}_{0}\left(  t\right)  +\mathcal{J}_{0}^{\ast}\left(  t\right)   &
\geq\mathcal{J}_{0}\left(  t\right)  \left(  1-a_{i}\right)  \geq0.
\end{align*}
They yield :%
\[
\left\vert I_{i}\left(  t\right)  -\frac{\mathcal{J}_{2}\left(  t\right)
}{\mathcal{J}_{0}\left(  t\right)  }\right\vert \leq\frac{2a_{i}%
\mathcal{J}_{2}\left(  t\right)  }{\left(  1-a_{i}\right)  \mathcal{J}%
_{0}\left(  t\right)  }%
\]
Turning back to (\ref{Sarko}) we get :%
\begin{align*}
&  \left\vert \mathbb{E}\left[  K\left(  \frac{\left\Vert X\right\Vert }%
{h}\right)  \left\langle X,e_{i}\right\rangle ^{2}\right]  -v\left(  h\right)
\right\vert \\
&  \leq\frac{2a_{i}}{1-a_{i}}\int_{0}^{h}t^{2}K\left(  \frac{t}{h}\right)
\frac{\mathcal{J}_{2}\left(  t\right)  }{\mathcal{J}_{0}\left(  t\right)
}f_{\left\Vert X\right\Vert }\left(  t\right)  dt+\int_{0}^{h}t^{2}K\left(
\frac{t}{h}\right)  \left(  \frac{\mathcal{J}_{2}\left(  t\right)
}{\mathcal{J}_{0}\left(  t\right)  }-\frac{\rho\left(  t\right)  }{t}\right)
f_{\left\Vert X\right\Vert }\left(  t\right)  dt
\end{align*}
Remind that, in order to alleviate notations we remove the index $i$
$\mathcal{J}_{0}\left(  t\right)  $ and $\mathcal{J}_{2}\left(  t\right)  .$
Assume that%
\begin{equation}
\limsup_{t\rightarrow0}\sup_{i}\frac{\mathcal{J}_{2}\left(  t\right)
}{\mathcal{J}_{0}\left(  t\right)  }\frac{t}{\rho\left(  t\right)  }\leq M
\label{palavas}%
\end{equation}
then we get%
\begin{align*}
\lim_{h\rightarrow0}\sup_{i}\left\vert \mathbb{E}\left[  K\left(
\frac{\left\Vert X\right\Vert }{h}\right)  \left\langle X,e_{i}\right\rangle
^{2}\right]  -v\left(  h\right)  \right\vert  &  \leq\left(  \frac{2a_{i}%
M}{1-a_{i}}+M+1\right)  \int_{0}^{h}t\rho\left(  t\right)  K\left(  \frac
{t}{h}\right)  f_{\left\Vert X\right\Vert }\left(  t\right)  dt\\
&  \leq M^{\prime}v\left(  h\right)
\end{align*}
where $M^{\prime}$ is some constant which does not depend on $i$ or $h$. This
finally entails Lemma \ref{chocolate}.In order to finish the proof we prove
(\ref{palavas}) now as a Lemma.

\begin{lemma}
\label{flo}We have%
\[
\limsup_{t\rightarrow0}\sup_{i}\frac{\mathcal{J}_{2}\left(  t\right)
}{\mathcal{J}_{0}\left(  t\right)  }\frac{t}{\rho\left(  t\right)  }\leq M
\]

\end{lemma}

\begin{proof}
We need to turn back to the proof of Lemma (\ref{prel}) from which we pick :%
\[
\frac{\mathcal{J}_{2}\left(  t\right)  }{\mathcal{J}_{0}\left(  t\right)
}=2\frac{\rho\left(  t\right)  }{t}\frac{\int_{0}^{t/\rho\left(  t\right)
}\sqrt{x}f_{\neq i}\left(  t-\rho\left(  t\right)  x\right)  dx}{\int
_{0}^{t/\rho\left(  t\right)  }\frac{1}{\sqrt{x}}f_{\neq i}\left(
t-\rho\left(  t\right)  x\right)  dx}%
\]
\bigskip And it suffices to get :%
\begin{equation}
\limsup_{t\rightarrow0}\sup_{i}\frac{\int_{0}^{t/\rho\left(  t\right)  }%
\sqrt{x}f_{\neq i}\left(  t-\rho\left(  t\right)  x\right)  dx}{\int
_{0}^{t/\rho\left(  t\right)  }\frac{1}{\sqrt{x}}f_{\neq i}\left(
t-\rho\left(  t\right)  x\right)  dx}\leq M \label{heaven}%
\end{equation}
We start with noting that%
\begin{align*}
&  \frac{\int_{0}^{t/\rho\left(  t\right)  }\sqrt{x}f_{\neq i}\left(
t-\rho\left(  t\right)  x\right)  dx}{\int_{0}^{t/\rho\left(  t\right)  }%
\frac{1}{\sqrt{x}}f_{\neq i}\left(  t-\rho\left(  t\right)  x\right)  dx}\\
&  =\frac{\int_{0}^{1}\sqrt{x}f_{\neq i}\left(  t-\rho\left(  t\right)
x\right)  dx+\int_{1}^{t/\rho\left(  t\right)  }\sqrt{x}f_{\neq i}\left(
t-\rho\left(  t\right)  x\right)  dx}{\int_{0}^{1}\frac{1}{\sqrt{x}}f_{\neq
i}\left(  t-\rho\left(  t\right)  x\right)  dx+\int_{1}^{t/\rho\left(
t\right)  }\frac{1}{\sqrt{x}}f_{\neq i}\left(  t-\rho\left(  t\right)
x\right)  dx}\\
&  \leq\frac{\int_{0}^{1}\frac{1}{\sqrt{x}}f_{\neq i}\left(  t-\rho\left(
t\right)  x\right)  dx+\int_{1}^{t/\rho\left(  t\right)  }\sqrt{x}f_{\neq
i}\left(  t-\rho\left(  t\right)  x\right)  dx}{\int_{0}^{1}\frac{1}{\sqrt{x}%
}f_{\neq i}\left(  t-\rho\left(  t\right)  x\right)  dx+\int_{1}%
^{t/\rho\left(  t\right)  }\frac{1}{\sqrt{x}}f_{\neq i}\left(  t-\rho\left(
t\right)  x\right)  dx}\\
&  \leq1+\frac{\int_{1}^{t/\rho\left(  t\right)  }\sqrt{x}f_{\neq i}\left(
t-\rho\left(  t\right)  x\right)  dx}{\int_{0}^{1}\frac{1}{\sqrt{x}}f_{\neq
i}\left(  t-\rho\left(  t\right)  x\right)  dx}.
\end{align*}
We deal with the denumerator%
\[
\int_{0}^{1}\frac{1}{\sqrt{x}}f_{\neq i}\left(  t-\rho\left(  t\right)
x\right)  dx\geq f_{\neq i}\left(  t-\rho\left(  t\right)  \right)  \int
_{0}^{1}\frac{1}{\sqrt{x}}dx=2f_{\neq i}\left(  t-\rho\left(  t\right)
\right)
\]
and%
\[
\frac{\int_{1}^{t/\rho\left(  t\right)  }\sqrt{x}f_{\neq i}\left(
t-\rho\left(  t\right)  x\right)  dx}{\int_{0}^{1}\frac{1}{\sqrt{x}}f_{\neq
i}\left(  t-\rho\left(  t\right)  x\right)  dx}\leq\frac{1}{2}\int_{1}%
^{t/\rho\left(  t\right)  }\sqrt{x}\frac{f_{\neq i}\left(  t-\rho\left(
t\right)  x\right)  }{f_{\neq i}\left(  t-\rho\left(  t\right)  \right)  }dx
\]
Then we we rewrite%
\begin{align*}
\frac{f_{\neq i}\left(  t-\rho\left(  t\right)  x\right)  }{f_{\neq i}\left(
t-\rho\left(  t\right)  \right)  }  &  =\frac{f_{\neq i}\left(  \left(
t-\rho\left(  t\right)  \right)  +\rho\left(  t-\rho\left(  t\right)  \right)
\frac{\rho\left(  t\right)  }{\rho\left(  t-\rho\left(  t\right)  \right)
}\left(  1-x\right)  \right)  }{f_{\neq i}\left(  t-\rho\left(  t\right)
\right)  }\\
&  =\frac{f_{\neq i}\left(  \left(  t-\rho\left(  t\right)  \right)
+\rho\left(  t-\rho\left(  t\right)  \right)  y_{t}\right)  }{f_{\neq
i}\left(  t-\rho\left(  t\right)  \right)  }%
\end{align*}
with $y_{t}=\frac{\rho\left(  t\right)  }{\rho\left(  t-\rho\left(  t\right)
\right)  }\left(  1-x\right)  .$ It is plain that, if $y_{t}=y$ does not
depend on $t$ that%
\[
\frac{f_{\neq i}\left(  \left(  t-\rho\left(  t\right)  \right)  +\rho\left(
t-\rho\left(  t\right)  \right)  y\right)  }{f_{\neq i}\left(  t-\rho\left(
t\right)  \right)  }\underset{t\rightarrow0}{\rightarrow}\exp(y)
\]
Since $\frac{\rho\left(  t\right)  }{\rho\left(  t-\rho\left(  t\right)
\right)  }\rightarrow1$ (see display 2.11.2 in Bingham, Goldie, Teugels
(1987)) and by Proposition 3.10.2 ibidem,%
\[
\frac{f_{\neq i}\left(  \left(  t-\rho\left(  t\right)  \right)  +\rho\left(
t-\rho\left(  t\right)  \right)  y_{t}\right)  }{f_{\neq i}\left(
t-\rho\left(  t\right)  \right)  }\underset{t\rightarrow0}{\rightarrow}%
\exp(1-x).
\]
From this remark, proving the Lemma finally comes down to proving%
\begin{equation}
\sup_{i,t}\int_{1}^{t/\rho\left(  t\right)  }\sqrt{x}\frac{f_{\neq i}\left(
t-\rho\left(  t\right)  x\right)  }{f_{\neq i}\left(  t\right)  \exp\left(
-x\right)  }\exp\left(  -x\right)  \exp dx\leq M. \label{thievery}%
\end{equation}
We focus on%
\[
\frac{f_{\neq i}\left(  t-\rho\left(  t\right)  x\right)  }{f_{\neq i}\left(
t\right)  \exp\left(  -x\right)  }.
\]
By the representation Theorem 3.10.8 in Bingham, Goldie, Teugels (1987) and
since all functions $f_{\neq i}$ have he same auxiliary function $\rho$ :%
\[
\frac{f_{\neq i}\left(  t-\rho\left(  t\right)  x\right)  }{f_{\neq i}\left(
t\right)  \exp\left(  -x\right)  }=\exp\left[  x-\int_{t-\rho\left(  t\right)
x}^{t}\frac{du}{\rho\left(  u\right)  }\right]  .
\]
Now it is easily seen that $\rho$ is continuous and nondecreasing in a
neighborhood of $0$ (see the remark about this fact within the proof of Lemma
\ref{prel}). Hence for $t$ small enough%
\[
-x\frac{\rho\left(  t\right)  }{\rho\left(  t-\rho\left(  t\right)  x\right)
}\leq-\int_{t-\rho\left(  t\right)  x}^{t}\frac{du}{\rho\left(  u\right)
}\leq-x
\]
and%
\[
\frac{f_{\neq i}\left(  t-\rho\left(  t\right)  x\right)  }{f_{\neq i}\left(
t\right)  \exp\left(  -x\right)  }\leq1
\]
for all $0\leq x\leq t/\rho\left(  t\right)  $ and for $t$ close to zero. At
last (\ref{thievery}) holds, hence (\ref{heaven}) hence Lemma \ref{flo}.
\end{proof}

\textbf{Proof of Theorem \ref{KT} :}

In order to finish the proof of the Theorem we have to cope now with
$\Gamma_{K}^{\#d}$ (see (\ref{di})) since Lemma \ref{chocolate} provides a
fair estimate with $\Gamma_{K}^{d}$. We have to prove that :%
\[
\left\Vert \Gamma_{K}^{\#d}\right\Vert _{\infty}=O\left(  v\left(  h\right)
\right)
\]

Actually we will prove that this bound is true in Hilbert-Schmidt norm since
the Hilbert-Schmidt norm of $\Gamma_{K}^{\#d}$ is easier to handle than its sup-norm.

We have (\ref{bo}) in mind and we start from :%
\[
J\left(  v\right)  =\frac{\int\int\frac{ut}{\sqrt{v^{2}-u^{2}-t^{2}}}%
f_{i}\left(  t\right)  f_{j}\left(  u\right)  f_{\neq ij}\left(  \sqrt
{v^{2}-u^{2}-t^{2}}\right)  1\hspace{-3pt}1_{\left\{  v\geq\sqrt{u^{2}+t^{2}%
}\right\}  }dudt}{\int\int\frac{v}{\sqrt{v^{2}-u^{2}-t^{2}}}f_{i}\left(
t\right)  f_{j}\left(  u\right)  f_{\neq ij}\left(  \sqrt{v^{2}-u^{2}-t^{2}%
}\right)  1\hspace{-3pt}1_{\left\{  v\geq\sqrt{u^{2}+t^{2}}\right\}  }dudt}%
\]
This display was obtained from (\ref{bo2}) by a change of variable.
Introducing again the function $\varphi_{i}$ (see above) :%
\[
f_{i}\left(  t\right)  f_{j}\left(  u\right)  =f_{i}\left(  0\right)
f_{j}\left(  0\right)  \left(  1+\varphi_{i}\left(  t\right)  \right)  \left(
1+\varphi_{j}\left(  u\right)  \right)
\]
hence%
\begin{align*}
J\left(  v\right)   &  =\frac{\int\int\frac{ut}{\sqrt{v^{2}-u^{2}-t^{2}}%
}\left(  1+\varphi_{i}\left(  t\right)  \right)  \left(  1+\varphi_{j}\left(
u\right)  \right)  f_{\neq ij}\left(  \sqrt{v^{2}-u^{2}-t^{2}}\right)
1\hspace{-3pt}1_{\left\{  v\geq\sqrt{u^{2}+t^{2}}\right\}  }dudt}{\int
\int\frac{v}{\sqrt{v^{2}-u^{2}-t^{2}}}\left(  1+\varphi_{i}\left(  t\right)
\right)  \left(  1+\varphi_{j}\left(  u\right)  \right)  f_{\neq ij}\left(
\sqrt{v^{2}-u^{2}-t^{2}}\right)  1\hspace{-3pt}1_{\left\{  v\geq\sqrt
{u^{2}+t^{2}}\right\}  }dudt}\\
&  =\frac{\int\int\frac{ut}{\sqrt{v^{2}-u^{2}-t^{2}}}\varphi_{i}\left(
t\right)  \varphi_{j}\left(  u\right)  f_{\neq ij}\left(  \sqrt{v^{2}%
-u^{2}-t^{2}}\right)  1\hspace{-3pt}1_{\left\{  v\geq\sqrt{u^{2}+t^{2}%
}\right\}  }dudt}{v\int\int\frac{\left(  1+\varphi_{i}\left(  t\right)
+\varphi_{j}\left(  u\right)  +\varphi_{i}\left(  t\right)  \varphi_{j}\left(
u\right)  \right)  }{\sqrt{v^{2}-u^{2}-t^{2}}}f_{\neq ij}\left(  \sqrt
{v^{2}-u^{2}-t^{2}}\right)  1\hspace{-3pt}1_{\left\{  v\geq\sqrt{u^{2}+t^{2}%
}\right\}  }dudt}%
\end{align*}
We see that%
\[
J\left(  v\right)  =\frac{\int\int\frac{ut}{\sqrt{v^{2}-u^{2}-t^{2}}}%
\varphi_{i}\left(  t\right)  \varphi_{j}\left(  u\right)  f_{\neq ij}\left(
\sqrt{v^{2}-u^{2}-t^{2}}\right)  1\hspace{-3pt}1_{\left\{  v\geq\sqrt
{u^{2}+t^{2}}\right\}  }dudt}{\int\int\frac{v}{\sqrt{v^{2}-u^{2}-t^{2}}%
}\left(  1+\varphi_{i}\left(  t\right)  \right)  \left(  1+\varphi_{j}\left(
u\right)  \right)  f_{\neq ij}\left(  \sqrt{v^{2}-u^{2}-t^{2}}\right)
1\hspace{-3pt}1_{\left\{  v\geq\sqrt{u^{2}+t^{2}}\right\}  }dudt}%
\]
since obviously%
\begin{align*}
\int\int\frac{ut}{\sqrt{v^{2}-u^{2}-t^{2}}}f_{\neq ij}\left(  \sqrt
{v^{2}-u^{2}-t^{2}}\right)  1\hspace{-3pt}1_{\left\{  v\geq\sqrt{u^{2}+t^{2}%
}\right\}  }dudt  &  =0\\
\int\int\frac{ut}{\sqrt{v^{2}-u^{2}-t^{2}}}\varphi_{i}\left(  t\right)
f_{\neq ij}\left(  \sqrt{v^{2}-u^{2}-t^{2}}\right)  1\hspace{-3pt}1_{\left\{
v\geq\sqrt{u^{2}+t^{2}}\right\}  }dudt  &  =0.
\end{align*}
We treat the numerator and the denumerator separatedly. Let%
\begin{align*}
\mathcal{N}  &  =\int\int\frac{ut}{\sqrt{v^{2}-u^{2}-t^{2}}}\varphi_{i}\left(
t\right)  \varphi_{j}\left(  u\right)  f_{\neq ij}\left(  \sqrt{v^{2}%
-u^{2}-t^{2}}\right)  1\hspace{-3pt}1_{\left\{  v\geq\sqrt{u^{2}+t^{2}%
}\right\}  }dudt\\
\mathcal{D}  &  =\int\int v\frac{\left(  1+\varphi_{i}\left(  t\right)
+\varphi_{j}\left(  u\right)  +\varphi_{i}\left(  t\right)  \varphi_{j}\left(
u\right)  \right)  }{\sqrt{v^{2}-u^{2}-t^{2}}}f_{\neq ij}\left(  \sqrt
{v^{2}-u^{2}-t^{2}}\right)  1\hspace{-3pt}1_{\left\{  v\geq\sqrt{u^{2}+t^{2}%
}\right\}  }dudt.
\end{align*}
Denoting again%
\[
a_{i}=\sup_{t\in\mathcal{V}_{0}}\left\vert \varphi_{i}\left(  t\right)
\right\vert ,
\]
we have
\[
\left\vert \mathcal{N}\right\vert \leq a_{i}a_{j}\int\int\frac{\left\vert
ut\right\vert }{\sqrt{v^{2}-u^{2}-t^{2}}}f_{\neq ij}\left(  \sqrt{v^{2}%
-u^{2}-t^{2}}\right)  1\hspace{-3pt}1_{\left\{  v\geq\sqrt{u^{2}+t^{2}%
}\right\}  }dudt
\]
and
\begin{align*}
\mathcal{D}  &  \geq\int\int v\frac{\left(  1-\left\vert \varphi_{i}\left(
t\right)  \right\vert -\left\vert \varphi_{j}\left(  u\right)  \right\vert
+\left\vert \varphi_{i}\left(  t\right)  \varphi_{j}\left(  u\right)
\right\vert \right)  }{\sqrt{v^{2}-u^{2}-t^{2}}}f_{\neq ij}\left(  \sqrt
{v^{2}-u^{2}-t^{2}}\right)  1\hspace{-3pt}1_{\left\{  v\geq\sqrt{u^{2}+t^{2}%
}\right\}  }dudt\\
&  \geq\int\int v\frac{\left(  1-a_{i}-a_{j}-a_{i}a_{j}\right)  }{\sqrt
{v^{2}-u^{2}-t^{2}}}f_{\neq ij}\left(  \sqrt{v^{2}-u^{2}-t^{2}}\right)
1\hspace{-3pt}1_{\left\{  v\geq\sqrt{u^{2}+t^{2}}\right\}  }dudt
\end{align*}
It is plain that for sufficiently large $i$ and $j,$%
\[
1-a_{i}-a_{j}-a_{i}a_{j}>0
\]
hence%
\[
J\left(  v\right)  \leq\frac{a_{i}a_{j}}{\left(  1-a_{i}-a_{j}-a_{i}%
a_{j}\right)  }\frac{\int\int\frac{\left\vert ut\right\vert }{\sqrt
{v^{2}-u^{2}-t^{2}}}f_{\neq ij}\left(  \sqrt{v^{2}-u^{2}-t^{2}}\right)
1\hspace{-3pt}1_{\left\{  v\geq\sqrt{u^{2}+t^{2}}\right\}  }dudt}{\int
\int\frac{v}{\sqrt{v^{2}-u^{2}-t^{2}}}f_{\neq ij}\left(  \sqrt{v^{2}%
-u^{2}-t^{2}}\right)  1\hspace{-3pt}1_{\left\{  v\geq\sqrt{u^{2}+t^{2}%
}\right\}  }dudt}%
\]

and $\lim_{i,j\rightarrow+\infty}a_{i}+a_{j}+a_{i}a_{j}=0$ hence is smaller
than $0.5$ for large $i$ and $j.$ Then we apply the same change of variable as
in (\ref{th}). We get for $i$ and $j$ large enough :%
\begin{align*}
J\left(  v\right)   &  \leq2a_{i}a_{j}\frac{\int_{0}^{v}\int_{0}^{2\pi}%
\frac{r^{3}\left\vert \cos\theta\sin\theta\right\vert }{\sqrt{v^{2}-r^{2}}%
}f_{\neq ij}\left(  \sqrt{v^{2}-r^{2}}\right)  drd\theta}{\int_{0}^{v}\int
_{0}^{2\pi}\frac{vr}{\sqrt{v^{2}-r^{2}}}f_{\neq ij}\left(  \sqrt{v^{2}-r^{2}%
}\right)  drd\theta}\\
&  \leq4\pi a_{i}a_{j}v\frac{\int_{0}^{1}\frac{x^{3}}{\sqrt{1-x^{2}}}f_{\neq
ij}\left(  v\sqrt{1-x^{2}}\right)  dx}{\int_{0}^{1}\frac{x}{\sqrt{1-x^{2}}%
}f_{\neq ij}\left(  v\sqrt{1-x^{2}}\right)  dx}%
\end{align*}
Now we invoke Lemma \ref{prel} to get :%
\[
\frac{\int_{0}^{1}\frac{x^{3}}{\sqrt{1-x^{2}}}f_{\neq ij}\left(
v\sqrt{1-x^{2}}\right)  dx}{\int_{0}^{1}\frac{x}{\sqrt{1-x^{2}}}f_{\neq
ij}\left(  v\sqrt{1-x^{2}}\right)  dx}\sim2\frac{\rho\left(  v\right)  }{v}%
\]
We should check now that :%
\[
\lim\sup_{v\rightarrow0}\sup_{i,j}\frac{v}{\rho\left(  v\right)  }\frac
{\int_{0}^{1}\frac{x^{3}}{\sqrt{1-x^{2}}}f_{\neq ij}\left(  v\sqrt{1-x^{2}%
}\right)  dx}{\int_{0}^{1}\frac{x}{\sqrt{1-x^{2}}}f_{\neq ij}\left(
v\sqrt{1-x^{2}}\right)  dx}<M
\]
The detailks of this step are omiied since they copy almost verbatim those of
Lemma \ref{flo}.

At last :%
\[
J\left(  v\right)  \leq10\pi a_{i}a_{j}\rho\left(  v\right)
\]
for large $i$ and $j$ and small $v$.

With this inequality in hand, we go back to (\ref{bo})%
\[
\mathbb{E}\left[  K\left(  \frac{\left\Vert X\right\Vert }{h}\right)
\left\langle X,e_{i}\right\rangle \left\langle X,e_{j}\right\rangle \right]
\leq10\pi a_{i}a_{j}\int_{0}^{h}v^{2}\rho\left(  v\right)  K\left(  \frac
{v}{h}\right)  f_{\left\Vert X\right\Vert }\left(  v\right)  dv
\]
from which we deduce that%
\[
\left\Vert \Gamma_{K}^{\#d}\right\Vert _{\infty}^{2}\leq\left\Vert \Gamma
_{K}^{\#d}\right\Vert _{2}^{2}\leq\left(  10\pi\int_{0}^{h}v^{2}\rho\left(
v\right)  K\left(  \frac{v}{h}\right)  f_{\left\Vert X\right\Vert }\left(
v\right)  dv\right)  ^{2}\left(  \sum_{i=1}^{+\infty}a_{i}^{2}\right)  ^{2}%
\]
and%
\[
\left\Vert \Gamma_{K}^{\#d}\right\Vert _{\infty}\leq10\pi\left(  \sum
_{i=1}^{+\infty}a_{i}^{2}\right)  \mathbb{E}\left(  \left\Vert X\right\Vert
^{2}\rho\left(  \left\Vert X\right\Vert \right)  K\left(  \frac{\left\Vert
X\right\Vert }{h}\right)  \right)  .
\]
Since $K$ has compact support we may say that $\left\Vert X\right\Vert \leq h$
hence :%
\begin{align*}
\mathbb{E}\left(  \left\Vert X\right\Vert ^{2}\rho\left(  \left\Vert
X\right\Vert \right)  K\left(  \frac{\left\Vert X\right\Vert }{h}\right)
\right)   &  \leq hv\left(  h\right)  ,\\
\left\Vert \Gamma_{K}^{\#d}\right\Vert _{\infty}  &  \leq10\pi\left(
\sum_{i=1}^{+\infty}a_{i}^{2}\right)  hv\left(  h\right)  .
\end{align*}
Together with Lemma \ref{chocolate} and assumption (\ref{mig}) this last
inequality yields Theorem \ref{KT}.

We are going to prove Proposition \ref{Compay} but first we give a Lemma :

\begin{lemma}
\label{free}Let $m,$ $p\in\mathbb{N}$ :%
\[
\mathbb{E}\left[  K^{m}\left(  \frac{\left\Vert X_{1}-x_{0}\right\Vert }%
{h}\right)  \left\Vert X_{1}-x_{0}\right\Vert ^{p}\right]  \sim K^{m}\left(
1\right)  h^{p}F\left(  h\right)
\]

\end{lemma}

\begin{proof}%
\begin{align*}
\mathbb{E}\left[  K^{m}\left(  \frac{\left\Vert X_{1}-x_{0}\right\Vert }%
{h}\right)  \left\Vert X_{1}-x_{0}\right\Vert ^{p}\right]   &  =\int_{0}%
^{h}u^{p}K^{m}\left(  \frac{u}{h}\right)  d\mathbb{P}^{\left\Vert X_{1}%
-x_{0}\right\Vert }\left(  u\right) \\
&  =h^{p}\int_{0}^{1}u^{p}K^{m}\left(  u\right)  d\mathbb{P}^{\left\Vert
X_{1}-x_{0}\right\Vert /h}\left(  u\right)
\end{align*}
We apply Fubini's theorem. It is plain that :%
\[
K^{m}\left(  u\right)  u^{p}=K^{m}\left(  1\right)  -\int_{u}^{1}\left[
s^{p}K^{m}\left(  s\right)  \right]  ^{\prime}ds,
\]
hence%
\begin{align*}
&  \int_{0}^{1}K^{m}\left(  u\right)  u^{p}d\mathbb{P}^{\left\Vert X_{1}%
-x_{0}\right\Vert /h}\left(  u\right) \\
&  =K^{m}\left(  1\right)  \int_{0}^{1}d\mathbb{P}^{\left\Vert X_{1}%
-x_{0}\right\Vert /h}\left(  u\right)  -\int\int\left[  K^{m}\left(  s\right)
s^{p}\right]  ^{\prime}1\hspace{-3pt}1_{\left\{  0\leq u\leq s\leq1\right\}
}d\mathbb{P}^{\left\Vert X_{1}-x_{0}\right\Vert /h}\left(  u\right)  ds\\
&  =K^{m}\left(  1\right)  F\left(  h\right)  -\int_{0}^{1}\left[  s^{p}%
K^{m}\left(  s\right)  \right]  ^{\prime}F\left(  hs\right)  ds\\
&  =K^{m}\left(  1\right)  F\left(  h\right)  \left(  1-\int_{0}^{1}\left[
s^{p}K^{m}\left(  s\right)  \right]  ^{\prime}\frac{F\left(  hs\right)
}{F\left(  h\right)  }ds\right)
\end{align*}
which finally gives%
\[
\mathbb{E}\left[  K^{m}\left(  \frac{\left\Vert X_{1}-x_{0}\right\Vert }%
{h}\right)  \left\Vert X_{1}-x_{0}\right\Vert ^{p}\right]  =h^{p}K^{m}\left(
1\right)  F\left(  h\right)  \left(  1-\int_{0}^{1}\left[  s^{p}K^{m}\left(
s\right)  \right]  ^{\prime}\frac{F\left(  hs\right)  }{F\left(  h\right)
}ds\right)  .
\]
We deal with $\int_{0}^{1}\left[  K^{m}\left(  s\right)  s^{p}\right]
^{\prime}\frac{F\left(  hs\right)  }{F\left(  h\right)  }ds$ and just have to
show that this integral goes to zero when $h$ does to prove the Lemma. Remind
that assumption $\mathbf{A}_{2}$ ensures that $K\left(  1\right)  >0$ and that
$\sup_{s}\left\vert K\left(  s\right)  ^{\prime}\right\vert <+\infty.$ Hence :%
\[
\sup_{s}\left\vert s^{p}K^{m}\left(  s\right)  ^{\prime}\right\vert <+\infty.
\]
At last Fact 1 (see display (\ref{F1})) ensures that $F\left(  hs\right)
/F\left(  h\right)  \rightarrow0$ when $s$ is fixed. Noticing that $F\left(
hs\right)  /F\left(  h\right)  \leq1$ and Lebesgue's dominated convergence
Theorem are enough to get the desired result and to complete the proof of the
Lemma.\bigskip
\end{proof}

\textbf{Proof of Proposition \ref{Compay} :} We begin with the first part of
the Proposition. We must prove that :%
\begin{align*}
v\left(  h\right)   &  =\mathbb{E}\left[  K\left(  \frac{\left\Vert
X_{1}-x_{0}\right\Vert }{h}\right)  \left\Vert X_{1}-x_{0}\right\Vert
\rho\left(  \left\Vert X_{1}-x_{0}\right\Vert \right)  \right] \\
&  =o\left(  \mathbb{E}\left[  K\left(  \frac{\left\Vert X_{1}-x_{0}%
\right\Vert }{h}\right)  \left\Vert X_{1}-x_{0}\right\Vert ^{2}\right]
\right)
\end{align*}

By (\ref{F2}) and as $\rho$ is a positive function, we can assume that $\rho$
is non-decreasing in a neighborhood of $0$ hence that :%
\begin{align*}
&  \mathbb{E}\left[  K\left(  \frac{\left\Vert X_{1}-x_{0}\right\Vert }%
{h}\right)  \left\Vert X_{1}-x_{0}\right\Vert \rho\left(  \left\Vert
X_{1}-x_{0}\right\Vert \right)  \right] \\
&  \leq\rho\left(  h\right)  \mathbb{E}\left[  K\left(  \frac{\left\Vert
X_{1}-x_{0}\right\Vert }{h}\right)  \left\Vert X_{1}-x_{0}\right\Vert \right]
\end{align*}
since the support of $K$ is $\left[  0,1\right]  .$ Then applying Lemma
\ref{free} with $m=1$ and $p=1$ we get%
\[
\mathbb{E}\left[  K\left(  \frac{\left\Vert X_{1}-x_{0}\right\Vert }%
{h}\right)  \left\Vert X_{1}-x_{0}\right\Vert \right]  \sim K\left(  1\right)
hF\left(  h\right)  .
\]
\bigskip Then for $h$ small enough :%
\begin{equation}
\frac{\mathbb{E}\left[  K\left(  \frac{\left\Vert X_{1}-x_{0}\right\Vert }%
{h}\right)  \left\Vert X_{1}-x_{0}\right\Vert \rho\left(  \left\Vert
X_{1}-x_{0}\right\Vert \right)  \right]  }{\mathbb{E}\left[  K\left(
\frac{\left\Vert X_{1}-x_{0}\right\Vert }{h}\right)  \left\Vert X_{1}%
-x_{0}\right\Vert ^{2}\right]  }\leq\frac{2K\left(  1\right)  \rho\left(
h\right)  hF\left(  h\right)  }{\mathbb{E}\left[  K\left(  \frac{\left\Vert
X_{1}-x_{0}\right\Vert }{h}\right)  \left\Vert X_{1}-x_{0}\right\Vert
^{2}\right]  } \label{vigne}%
\end{equation}
By Lemma \ref{free} again with $m=1$ and $p=2,$ the denumerator behaves like
$K\left(  1\right)  h^{2}F\left(  h\right)  .$ At last since $\rho\left(
h\right)  /h\rightarrow0,$ the display above tends to zero.\newline We go on
with the second part of the Proposition. It is not difficult by copying the
proof of Lemma \ref{free} to show that :%
\begin{equation}
v\left(  h\right)  =K\left(  1\right)  \rho\left(  h\right)  hF\left(
h\right)  \left(  1-\int_{0}^{1}\frac{\left[  sK\left(  s\right)  \rho\left(
hs\right)  \right]  ^{\prime}}{\rho\left(  h\right)  }\frac{F\left(
hs\right)  }{F\left(  h\right)  }ds\right)  \label{toto}%
\end{equation}
and%
\[
\frac{\left[  sK\left(  s\right)  \rho\left(  hs\right)  \right]  ^{\prime}%
}{\rho\left(  h\right)  }=\frac{K\left(  s\right)  \rho\left(  hs\right)
}{\rho\left(  h\right)  }+\frac{sK^{\prime}\left(  s\right)  \rho\left(
hs\right)  }{\rho\left(  h\right)  }+\frac{hsK\left(  s\right)  \rho^{\prime
}\left(  hs\right)  }{\rho\left(  h\right)  }.
\]
The first and second term on the right are uniformly bounded with respect to
$h$ and $s$ since $\rho$ is non decreasing in a neighborhood of zero and $K$
and $K^{\prime}$ are bounded. We turn to the last. Remind that here $\rho$ is
assumed to be regularly varying at zero. From (\ref{F2}) we now that the index
of regular variation of $\rho$ is $d\geq0.$ Now open the book by Bingham,
Goldie and Teugels (1987). The definition of regular variation is given p.18.
from Theorem 1.7.2 p.39 we deduce that $\rho^{\prime}$ is regularly varying
with index $d-1$ and also that :%
\[
\limsup_{t\rightarrow0}\left\vert \frac{t\rho^{\prime}\left(  t\right)  }%
{\rho\left(  t\right)  }\right\vert =d
\]
which means that $hsK\left(  s\right)  \rho^{\prime}\left(  hs\right)
/\rho\left(  h\right)  $ is unformly bounded with respect to $h$ and $s$ for
small $h$ and $0\leq s\leq1$. Applying the dominated convergence theorem to
(\ref{toto}) yields the desired result.\bigskip

\textbf{Proof of Theorem \ref{covop} :}

Simple calculations give :%
\begin{align*}
&  \mathbb{E}\left\Vert \Gamma_{K,n}-\Gamma_{K}\right\Vert ^{2}\\
&  =\frac{1}{n}\left(  \mathbb{E}\left[  K^{2}\left(  \frac{\left\Vert
X_{1}-x_{0}\right\Vert }{h}\right)  \left\Vert X_{1}-x_{0}\right\Vert
^{4}\right]  -\left\Vert \Gamma_{K}\right\Vert ^{2}\right)
\end{align*}

We begin with computing $\mathbb{E}\left[  K^{2}\left(  \frac{\left\Vert
X_{1}-x_{0}\right\Vert }{h}\right)  \left\Vert X_{1}-x_{0}\right\Vert
^{4}\right]  .$ Lemma \ref{free} gives :%
\[
\mathbb{E}\left[  K^{2}\left(  \frac{\left\Vert X_{1}-x_{0}\right\Vert }%
{h}\right)  \left\Vert X_{1}-x_{0}\right\Vert ^{4}\right]  \sim K^{2}\left(
1\right)  h^{4}F\left(  h\right)
\]

Under the assumptions of Theorem \ref{KT} we also have :%
\begin{equation}
\left\Vert \Gamma_{K}\right\Vert =O\left(  v\left(  h\right)  \right)
\end{equation}
where (see (\ref{vigne})) :%
\[
\frac{v\left(  h\right)  }{K\left(  1\right)  \rho\left(  h\right)  hF\left(
h\right)  }\rightarrow1
\]
and%
\[
\left\Vert \Gamma_{K}\right\Vert ^{2}=O\left(  \left[  \rho\left(  h\right)
hF\left(  h\right)  \right]  ^{2}\right)  =o\left(  h^{4}K^{2}\left(
1\right)  F\left(  h\right)  \right)
\]
which also means that%
\begin{align*}
\mathbb{E}\left\Vert \Gamma_{K,n}-\Gamma_{K}\right\Vert ^{2}  &  \sim\frac
{1}{n}\mathbb{E}\left[  K^{2}\left(  \frac{\left\Vert X_{1}-x_{0}\right\Vert
}{h}\right)  \left\Vert X_{1}-x_{0}\right\Vert ^{4}\right] \\
&  \sim K^{2}\left(  1\right)  \frac{h^{4}F\left(  h\right)  }{n}.
\end{align*}

Theorem \ref{covop} is now proved.\bigskip

\textbf{Proof of Corollary \ref{spec} :}

The proof of the first display, related with the eigenvalues, stems from the
famous bound :%
\[
\sup_{p\in\mathbb{N}}\left\vert s_{p}\left(  T\right)  -s_{p}\left(  U\right)
\right\vert \leq\left\Vert U-T\right\Vert _{\infty}%
\]
where $U$ and $Y$ are compact linear operators from and onto $H$ and
$s_{p}\left(  T\right)  $ stands for the $p^{th}$ characteristic number of the
operator $T$. The proof of the second display, namely on eigenprojections,
will be deduced from an article by Mas and Menneteau (2003) as announced
sooner in the paper.\bigskip

\textbf{Proof of Proposition \ref{Juliet} : }We want to estimate $P\left(
\left\Vert X\right\Vert \leq\varepsilon\right)  $ when $\varepsilon$ goes to
$0.$ In order to understand these few lines, the reader is referred to Dembo
et alii (1995). We have to compute or to give asymptotic equivalents for
formulas (3), (4) and (5) in their article. Let $\theta=\theta\left(
\varepsilon\right)  $ be the unique solution of :%
\[
\mu\left(  \theta\right)  =\varepsilon
\]
where here%
\[
\mu\left(  \theta\right)  =\sum_{i=1}^{+\infty}\frac{1}{\exp\left(  \alpha
i\right)  +2\theta}.
\]
We begin with :%
\begin{align}
\mu\left(  \theta\right)   &  =\sum_{i=1}^{+\infty}\frac{1}{\exp\left(  \alpha
i\right)  +2\theta}\sim\int_{1}^{+\infty}\frac{dx}{\exp\left(  \alpha
x\right)  +2\theta}\nonumber\\
&  =\frac{1}{\alpha}\int_{0}^{\exp\left(  -\alpha\right)  }\frac{du}{1+2\theta
u}=\frac{\log\left(  1+2\theta\exp\left(  -\alpha\right)  \right)  }%
{2\alpha\theta}. \label{terra}%
\end{align}

Then%
\begin{align*}
\psi\left(  \theta\right)   &  =\sqrt{2\sum_{i=1}^{+\infty}\left(
\frac{\theta}{\exp\left(  \alpha i\right)  +2\theta}\right)  ^{2}}\sim
\theta\sqrt{2\int_{1}^{+\infty}\frac{dx}{\left(  \exp\left(  \alpha x\right)
+2\theta\right)  ^{2}}}\\
&  =\theta\sqrt{2\frac{1}{\alpha}\int_{0}^{\exp\left(  -\alpha\right)  }%
\frac{udu}{\left(  1+2\theta u\right)  ^{2}}}\\
&  =\sqrt{\frac{1}{2\alpha}}\sqrt{\log\left(  1+2\theta\exp\left(
-\alpha\right)  \right)  -\frac{2\theta\exp\left(  -\alpha\right)  }%
{1+2\theta\exp\left(  -\alpha\right)  }}%
\end{align*}

At last%
\begin{equation}
I\left(  \theta\right)  =\frac{1}{2}\sum_{i=1}^{+\infty}\log\left(
1+2\theta\exp\left(  -\alpha i\right)  \right)  -\theta\mu\left(
\theta\right)  . \label{ekova}%
\end{equation}
We focus on%
\[
\frac{1}{2}\sum_{i=1}^{+\infty}\log\left(  1+2\theta\exp\left(  -\alpha
i\right)  \right)  \sim\frac{1}{2}\int_{1}^{+\infty}\log\left(  1+2\theta
\exp\left(  -\alpha x\right)  \right)  dx.
\]
Setting $u=2\theta\exp\left(  -\alpha x\right)  $ this last integral becomes%
\begin{align*}
&  \frac{1}{2}\int_{1}^{+\infty}\log\left(  1+2\theta\exp\left(  -\alpha
x\right)  \right)  dx\\
&  =\frac{1}{2\alpha}\int_{0}^{2\theta\exp\left(  -\alpha\right)  }\frac
{\log\left(  1+u\right)  }{u}du\\
&  \sim\frac{1}{2\alpha}\int_{1}^{2\theta\exp\left(  -\alpha\right)  }%
\frac{\log\left(  1+u\right)  }{u}du\\
&  \sim\frac{1}{2\alpha}\int_{1}^{2\theta\exp\left(  -\alpha\right)  }%
\frac{\log\left(  u\right)  }{u}du=\frac{1}{4\alpha}\left[  \log\left(
2\theta\exp\left(  -\alpha\right)  \right)  \right]  ^{2}\\
&  \sim\frac{1}{4\alpha}\left[  \log\theta\right]  ^{2}.
\end{align*}
From (\ref{ekova}) we see that :%
\[
\frac{I\left(  \theta\right)  }{\frac{1}{4\alpha}\left[  \log\theta\right]
^{2}}=p\left(  \theta\right)  -\frac{\theta\mu\left(  \theta\right)  }%
{\frac{1}{4\alpha}\left[  \log\theta\right]  ^{2}},
\]
where $p\left(  \theta\right)  $ tends to $1$ when $\theta$ goes to infinity.
Then from (\ref{terra}) we get :%
\[
\frac{\theta\mu\left(  \theta\right)  }{\frac{1}{4\alpha}\left[  \log
\theta\right]  ^{2}}\sim\frac{4\alpha}{\log\theta}\rightarrow0
\]
and
\[
\frac{I\left(  \theta\right)  }{\frac{1}{4\alpha}\left[  \log\theta\right]
^{2}}\rightarrow1.
\]

Collecting our results we have the following final asymptotic equivalence%
\begin{align*}
\theta &  \sim-\frac{\log\left(  \varepsilon\right)  }{2\alpha\varepsilon}\\
\mu\left(  \theta\right)   &  =\varepsilon\sim\frac{\log\left(  \theta\right)
}{2\alpha\theta}\\
\psi\left(  \theta\right)   &  \sim\sqrt{\frac{1}{2\alpha}\log\left(
\theta\right)  }\sim\sqrt{\frac{-\log\left(  \varepsilon\right)  }{2\alpha}}\\
I\left(  \theta\right)   &  \sim\frac{1}{4\alpha}\left[  \log\left(
\theta\right)  \right]  ^{2}\sim\frac{1}{4\alpha}\left[  \log\left(
\varepsilon\right)  \right]  ^{2}%
\end{align*}
Collecting these formulas together with display (10) in the above-mentioned
article we get :%
\[
\mathbb{P}\left(  \left\Vert X\right\Vert <\varepsilon\right)  \sim\sqrt
{\frac{\alpha}{-\pi\log\left(  \varepsilon\right)  }}\exp\left(  -\frac
{1}{4\alpha}\left[  \log\left(  \varepsilon\right)  \right]  ^{2}\right)  .
\]

\bigskip\bigskip\textbf{Acknowledgements} : I am sincerely grateful to S.
Ga\"{\i}ffas who pointed out to me the existence of de Haan's class $\Gamma
.$

\end{document}